\pdfoutput=1
\documentclass[12pt]{article}
\usepackage{amsmath,amssymb,amsthm}
\usepackage{mathtools}
\usepackage{array}
\usepackage{geometry}
\usepackage{xurl}
\usepackage{hyperref}
\hypersetup{
  pdftitle={A Lean 4 Framework for the Radii Polynomial Method},
  pdfauthor={Fengyang Wang, Konstantin Mischaikow, Marcio Gameiro},
  pdfsubject={Computer-assisted proofs; formal verification},
  pdfkeywords={radii polynomial method, computer-assisted proof, Lean 4, Mathlib, weighted sequence algebras, initial value problems},
  colorlinks=true, linkcolor=blue!60!black, citecolor=blue!60!black, urlcolor=blue!60!black
}
\usepackage[expansion=false]{microtype}
\usepackage{xcolor}
\usepackage[perpage]{footmisc}
\usepackage{tikz-cd}
\usepackage{listings}
\makeatletter
\newcommand{\Needspace}[1]{\par\penalty-100\begingroup\dimen@=#1\relax\vskip\dimen@\penalty\z@\vskip-\dimen@\kern\z@\endgroup}
\makeatother
\lstdefinelanguage{lean4}{morekeywords={def,theorem,structure,class,abbrev,where},sensitive=true,morecomment=[l]{--}}
\newtheorem{theorem}{Theorem}[section]

\newtheorem{definition}[theorem]{Definition}

\newtheorem*{theorem*}{Theorem}

\newcommand{\N}{\mathbb{N}}
\newcommand{\Z}{\mathbb{Z}}
\newcommand{\Q}{\mathbb{Q}}
\newcommand{\R}{\mathbb{R}}

\newcommand{\xbar}{\bar{x}}
\newcommand{\abar}{\bar{a}}
\newcommand{\xtilde}{\tilde{x}}
\newcommand{\atilde}{\tilde{a}}
\newcommand{\Ynorm}{Y_0}
\newcommand{\Znorm}[1]{Z_{#1}}
\newcommand{\closedball}[2]{\overline{B}_{#1}(#2)}
\newcommand{\lean}[1]{\texttt{#1}}

\newcommand{\MvPoly}{\lean{MvPolynomial}}

\newcommand{\Pow}{\mathrm{Pow}}

\newcommand{\LibCommit}{a926ce1525a1}

\title{A Lean~4 Framework for the Radii Polynomial Method}
\author{\begin{tabular}[t]{@{}c@{\hspace{1em}}c@{\hspace{1em}}c@{}}
  Fengyang Wang & Konstantin Mischaikow & Marcio Gameiro\\[4pt]
  \small Department of Mathematics & \small Department of Mathematics & \small Department of Mathematics\\
  \small Rutgers University & \small Rutgers University & \small Rutgers University\\
  \small Piscataway, NJ 08854, USA & \small Piscataway, NJ 08854, USA & \small Piscataway, NJ 08854, USA\\
  \small \texttt{fengyang.w@rutgers.edu} & \small \texttt{mischaik@math.rutgers.edu} & \small \texttt{gameiro@math.rutgers.edu}
\end{tabular}}
\date{}

\begin{document}
\maketitle

\begin{abstract}
Understanding nonlinear dynamics is central to many fields of study, yet
most nonlinear systems are beyond the reach of closed-form solutions.
Numerical simulation is then used for revealing dynamical
structures before their existence can be
established. This alone, however, leaves the validity of the computed
solution unanswered. Rigorous numerical methods address this gap by
combining computation with explicit analytic estimates. Among them, the
radii polynomial approach, a Newton--Kantorovich type argument, turns a
numerical approximation into a proof of existence and uniqueness of a true
solution near that approximation. In previous practice the verification of these
inequalities relies on interval-arithmetic libraries and hand-derived
algebraic bounds. In this paper, we present a framework for the radii
polynomial method in the Lean~4 proof assistant, which turns numerical
certificates into fully type-checked theorems about the true solution. A
central challenge in this process is to manage the complexity of the
formalization. Our approach begins with a family of problems and
their common structure. The spaces in which
the analysis takes place are developed through universal properties in
suitable categories, with extension properties that characterize maps by
their action on generators. These properties guide the design of a reusable
interface such that problem-specific
arguments concentrate on the numerical data and the algebraic bounds that
distinguish one equation from another. We formalize a general radii
polynomial theorem, develop weighted coefficient algebras and reusable operator estimates, connect finite
computational checks to the hypotheses of the theorem, and establish the
realization theorems that transport certified coefficient space zeros to
analytic solutions of the original equations. We further indicate how the universal construction extends to weighted \(\ell^p\) spaces.
\end{abstract}



\section{Introduction}\label{sec:intro}

In applications, nonlinear problems are typically solved numerically.
The obvious concern is whether the proposed solution is an adequate
approximation of a true solution as opposed to a spurious solution.
Addressing these concerns has led to the development of a wide variety of
techniques, often referred to as computer-assisted proofs, that take as
input purported numerical solutions and provide as output either a guarantee
that a true solution lies within an explicit distance of the numerical
solution, or a statement that such a guarantee could not be provided
(see~\cite{vdbLessard2015,gomezSerrano2019,Nakao2019Numerical} and
references therein). Computer-assisted proofs of this kind go back to the 1982 proof of the
Feigenbaum conjectures~\cite{lanford1982}. Mischaikow and Mrozek established
chaos in the Lorenz equations~\cite{mischaikowMrozek1995}, with Szymczak at
the classical parameter values~\cite{mischaikowMS2001}, and Tucker certified
the Lorenz attractor in his resolution of Smale's 14th
problem~\cite{tucker2002}. These proofs combine interval arithmetic with
analytic estimates. In our development, the analytic estimates and the
numerical inequalities are Lean theorems. Lean~4~\cite{deMouraUllrich2021}
is a proof assistant, a program that checks every step of a proof against a
small kernel, and Mathlib~\cite{mathlib2020} is its library of
formalized mathematics, on which our development is built. A theorem
accepted by Lean is checked down to the axioms it rests on.

The focus of this paper is on the radii polynomial
method~\cite{dayLM2007,castelliGL2017,vdbMJLM2011}, which provides a
systematic framework by which a
Newton--Kantorovich type argument ensures the existence of a solution with
explicit bounds. In this setting the
computer is used to check certain inequalities. If they are satisfied, then
the existence of the desired solution is guaranteed. The goal of this paper
is to initiate the mathematical formalization of these results using
Lean~4.

In short, the radii polynomial method is this. Let \(F\colon X\to Y\) be a differentiable map
from a Banach space to a normed space, let \(\xbar\in X\) be a numerical approximation of a
zero of \(F\), let \(A^\dagger\colon X\to Y\) approximate the derivative
\(DF(\xbar)\), and let \(A\colon Y\to X\) be an injective approximate inverse
of \(A^\dagger\). Fix a radius \(r_0>0\), and suppose that numbers
\(\Ynorm,\Znorm{0},\Znorm{1}\) and a function \(\Znorm{2}\) of the radius
satisfy
\[
 \|A(F(\xbar))\|\le\Ynorm,\qquad
 \|I-A\circ A^\dagger\|\le\Znorm{0},\qquad
 \|A\circ[A^\dagger-DF(\xbar)]\|\le\Znorm{1},
\]
\[
 \|A\circ[DF(c)-DF(\xbar)]\|\le\Znorm{2}(r_0)\,r_0
 \qquad\text{whenever }\|c-\xbar\|\le r_0 .
\]
If the \emph{radii polynomial}
\[
 p(r)=\Znorm{2}(r)\,r^2-(1-\Znorm{0}-\Znorm{1})\,r+\Ynorm
\]
is negative at \(r_0\), then \(F\) has exactly one zero in the closed
ball \(\closedball{r_0}{\xbar}\) (Theorem~\ref{thm:radii-poly}).

The main application of this method is to differential equations, formulated as
zero-finding problems. To put this work in context, we consider two simple but classical problems.
\begin{description}
  \item[P1] Let \(f\colon\R\to\R\) be given by a polynomial. Assume
  \(\xbar\in\R\) is a purported zero of \(f\), i.e., \(f(\xbar)\approx0\).
  Seek a radius \(r_0>0\) for which the four bounds certify a
  unique zero \(\xtilde\) in \(\closedball{r_0}{\xbar}\).
  \item[P2] Consider a numerically generated solution \(\xbar\) of an initial
  value problem for a system of ordinary differential equations with a
  polynomial nonlinearity. More precisely, assume \(\xbar\) approximates a
  solution of \(\dot x=f(x)\), \(x(0)=x_0\), in a chosen normed function
  space. Seek a radius \(r_0>0\) for which the four bounds certify
  a unique solution \(\xtilde\in\closedball{r_0}{\xbar}\) of this initial
  value problem.
\end{description}
The fundamental difference between \textbf{P1} and \textbf{P2} is that in
\textbf{P1} \(\xtilde\in\R\), whereas in \textbf{P2} \(\xtilde\) is an
element of some function space.

In the power series examples of this paper, \textbf{P2} among them,
the unknown \(x\) is represented by the sequence
\(a=(a_n)_{n\ge0}\) of its power series \(x(t)=\sum_{n\ge0}a_nt^n\), and the
equation becomes an equation \(F(a)=0\) for the sequence \(a\). This turns
\textbf{P2} into a problem of \textbf{P1} type, in loose language.
The unknown is now an infinite sequence rather than a single number.
This paper formalizes the same zero-finding argument
in Lean on the weighted sequence space \(\ell^1_\nu\), while carrying over
Mathlib's existing polynomial rings and algebra structures.

The organizing idea of the library is to construct maps out of the
coefficient space from their values on the atoms, the sequences with one
coefficient equal to \(1\) and all others zero. If the norm of each assigned value is at most a common constant times the
norm of its atom, the construction gives a bounded linear operator
(Section~\ref{sec:free-banach}). When the assigned values also preserve
the unit and products, the resulting map is an algebra homomorphism
(Section~\ref{sec:monoid-algebra}), as in the evaluation of series
(Section~\ref{sec:evaluations}). In particular, polynomial substitution
commutes with evaluation. For \(\nu>0\), the space \(\ell^1_\nu\) consists
of the sequences with \(\sum_n|a_n|\,\nu^n<\infty\). Writing
\(\operatorname{ev}_t(a)=\sum_na_nt^n\), for a polynomial \(P\) in one
variable and a sequence \(a\in\ell^1_\nu\),
\[
 \operatorname{ev}_t\bigl(P(a)\bigr)=P\bigl(\operatorname{ev}_t(a)\bigr)
 \qquad(|t|\le\nu).
\]
Here \(P(a)\) uses the Cauchy product of sequences.
Section~\ref{sec:free} gives the same identity for a polynomial in several
variables, with a coefficient sequence substituted for each variable. The overhead of one problem then
concentrates on the numerical data and on the algebraic bounds that
distinguish one equation from another.

A \emph{certificate} consists of finitely many rational numbers together with
the proofs needed to apply Theorem~\ref{thm:radii-poly}, including the
four bounds and \(p(r_0)<0\). For the families treated here, a
\emph{family theorem} uses the certificate to give a unique coefficient zero
in the ball. A \emph{realization} theorem then carries the coefficient zero
back to a solution of the original equation. For the Taylor initial value
problems, realization uses evaluation and termwise differentiation.

The paper is organized as follows. Section~\ref{sec:theorem} states and
proves the radii polynomial theorem for maps from a Banach space to a normed
space, and shows its Lean form.
Section~\ref{sec:algebras} presents the weighted coefficient algebra in
Lean, the maps out of it, and polynomial substitution into it.
Appendix~\ref{app:categorical} gives the categorical formulation of
these maps and their composition.
Section~\ref{sec:reductions} reduces the four norm bounds to finite
calculations on coefficients and polynomial expressions for the
polynomial-zero and Taylor families.
Section~\ref{sec:examples} runs the certificates for a square-root branch
given by a power series and two Taylor initial value problems, among them the Lorenz system. Each certificate gives a unique coefficient
zero in its certified ball. The realization theorems give the square-root
identity and, for the initial value problems, solutions unique among those
with values in the stated ball.
Section~\ref{sec:verification} explains how Lean checks the finite
calculations and which axioms the resulting theorems depend on.
Section~\ref{sec:discussion} discusses reuse and further applications,
including the case \(f(x)=x^2-2\) of \textbf{P1}.
Appendix~\ref{app:lp} treats convolution and operator bounds on
weighted \(\ell^p\) spaces.

Our Lean~4 development\footnote{Source code at \url{https://github.com/IlPreteRosso/RadiiPolynomial}, revision \lean{\LibCommit}, built with Lean~4.34.1, Mathlib \lean{d13f23b7} (tag \lean{v4.34.1}), and LeanCert~\cite{leancert} \lean{7f91b6eb}.}
takes the numerical approximations as external data, and the certificate
checks verify that they satisfy the hypotheses of the radii polynomial
theorem.

\section{Radii polynomial theorem}\label{sec:theorem}

Let \(\xbar\in X\) be a numerical approximation of a zero of a nonlinear map
\(F\colon X\to Y\). A certificate supplies a radius \(r_0>0\) and proves four norm bounds.
Under the differentiability and injectivity hypotheses below,
Theorem~\ref{thm:radii-poly} gives exactly one zero in
\(\closedball{r_0}{\xbar}\) when the radii polynomial is negative at \(r_0\).
Theorem~\ref{thm:radii-poly} is a standard form of the radii polynomial theorem (\cite[Prop.~2.4]{castelliGL2017}, \cite[Thm.~1.6]{lessardMJ2017}).

\begin{theorem}[Radii polynomial theorem, \lean{general\_radii\_polynomial\_theorem}]\label{thm:radii-poly}
Let \(X\) be a real Banach space and let \(\xbar\in X\). Let \(Y\) be a
real normed linear space.\footnote{The core theorems of the library write \(E\), \(F\), \(f\), \lean{A\_dagger}, \(Df\), \lean{xBar}, and \lean{xTilde} for \(X\), \(Y\), \(F\), \(A^\dagger\), \(DF\), \(\xbar\), and \(\xtilde\). This paper uses its own notation. Elsewhere \(E\) denotes a fiber family and \(f\) a vector field.}
Let \(F\colon X\to Y\) be a Fr\'echet differentiable function with
Fr\'echet derivative \(DF\). Let \(A^\dagger\colon X\to Y\) and
\(A\colon Y\to X\) be continuous linear maps, and assume \(A\) is
injective. Fix a real number \(r_0>0\), constants
\(\Ynorm,\Znorm{0},\Znorm{1}\in\R\), and a function
\(\Znorm{2}\colon\R\to\R\). If
\begin{align}
  \|A(F(\xbar))\| &\leq \Ynorm,\label{eq:Y}\\
  \|I_X-A\circ A^\dagger\| &\leq\Znorm{0},\label{eq:Z0}\\
  \|A\circ[A^\dagger-DF(\xbar)]\| &\leq\Znorm{1},\label{eq:Z1}\\
  \|A\circ[DF(c)-DF(\xbar)]\| &\leq\Znorm{2}(r_0)\,r_0
       \quad(\forall c\in\closedball{r_0}{\xbar}),\label{eq:Z2}
\end{align}
and
\[
  p(r_0)=\Znorm{2}(r_0)\,r_0^2-(1-\Znorm{0}-\Znorm{1})\,r_0+\Ynorm<0,
\]
then there exists a unique \(\xtilde\in\closedball{r_0}{\xbar}\) with
\(F(\xtilde)=0\).
\end{theorem}

The map \(A\) is the \emph{preconditioner}. Condition~\eqref{eq:Y} bounds the preconditioned residual \(A(F(\xbar))\)
by \(\Ynorm\). The second condition~\eqref{eq:Z0} indicates that \(A\) is an approximate inverse of \(A^\dagger\), in the one-sided sense that
\(A\circ A^\dagger\) is close to \(I_X\). The third
condition~\eqref{eq:Z1} bounds the error \(A^\dagger-DF(\xbar)\) of the
\emph{approximate derivative} \(A^\dagger\), after composition with \(A\). Finally, the fourth
condition~\eqref{eq:Z2} bounds the variation of \(A\circ DF\) over the ball
\(\closedball{r_0}{\xbar}\). The four bounds define the radii polynomial
\[
  p(r)=\Znorm{2}(r)\,r^2-(1-\Znorm{0}-\Znorm{1})\,r+\Ynorm ,
\]
and the proof studies the Newton-like map \(T(x):=x-A(F(x))\) on the ball
\(\closedball{r_0}{\xbar}\subset X\). When \(p(r_0)<0\), the four bounds make
\(T\) a contraction of the ball into itself, and the contraction mapping
theorem gives a unique fixed point \(\xtilde\in\closedball{r_0}{\xbar}\) of
\(T\). Injectivity of \(A\), which the four bounds do not imply, gives \(T(x)=x\iff F(x)=0\), hence a unique zero of
\(F\) in \(\closedball{r_0}{\xbar}\).

The certificates of Section~\ref{sec:examples} take \(\Znorm{2}\) to be constant as a function of the radius.

In Lean the theorem reads
\begin{lstlisting}
theorem general_radii_polynomial_theorem
  {f : E → F} {xBar : E} {A : F →L[ℝ] E} {A_dagger : E →L[ℝ] F}
  {Y₀ Z₀ Z₁ : ℝ} {Z₂ : ℝ → ℝ} {r₀ : ℝ}
  (hr₀ : 0 < r₀)
  (h_Y₀ : ‖A (f xBar)‖ ≤ Y₀)
  (h_Z₀ : ‖I_E - A.comp A_dagger‖ ≤ Z₀)
  (h_Z₁ : ‖A.comp (A_dagger - fderiv ℝ f xBar)‖ ≤ Z₁)
  (h_Z₂ : ∀ c ∈ Metric.closedBall xBar r₀,
    ‖A.comp (fderiv ℝ f c - fderiv ℝ f xBar)‖ ≤ Z₂ r₀ * r₀)
  (hf_diff : Differentiable ℝ f)
  (h_radii : generalRadiiPolynomial Y₀ Z₀ Z₁ Z₂ r₀ < 0)
  (hA_inj : Function.Injective A) :
  ∃! xTilde ∈ Metric.closedBall xBar r₀, f xTilde = 0
\end{lstlisting}
The operators \(A\) and \(A^\dagger\) of Theorem~\ref{thm:radii-poly} are
the arguments \lean{A} and \lean{A\_dagger}, bundled continuous linear
maps. The four hypotheses \lean{h\_Y\(_0\)} to \lean{h\_Z\(_2\)} are the four
conditions~\eqref{eq:Y}--\eqref{eq:Z2}, with Mathlib's Fr\'echet derivative
\lean{fderiv} in the role of \(DF\). The hypothesis \lean{h\_radii} is
the radii inequality \(p(r_0)<0\), where \lean{generalRadiiPolynomial} is the radii polynomial
\(p\), and the conclusion is the unique zero in the closed
ball. For one example, the certificate and the family theorem that consumes it together supply the hypotheses.

Theorem~\ref{thm:radii-poly} assumes differentiability on all of \(X\),
although its derivative estimates concern only the closed ball.

\begin{theorem}[Local form,
  \lean{general\_radii\_polynomial\_theorem\_local}]\label{thm:radii-local}
Keep the hypotheses of Theorem~\ref{thm:radii-poly}, except that \(F\) need
not be differentiable on \(X\). Instead, a continuous linear map
\(DF(x)\colon X\to Y\) is given for every \(x\in\closedball{r_0}{\xbar}\),
such that \(F\) has derivative \(DF(x)\) within \(\closedball{r_0}{\xbar}\)
at \(x\), and \eqref{eq:Y}--\eqref{eq:Z2} hold with this \(DF\). Then there exists a
unique \(\xtilde\in\closedball{r_0}{\xbar}\) with \(F(\xtilde)=0\).
\end{theorem}

Theorem~\ref{thm:radii-local} is the more general statement.
Theorem~\ref{thm:radii-poly} takes the derivative within the ball from
differentiability on \(X\) and is its direct consequence. Under the same hypotheses, the zero lies in the
open ball and \(DF(\xtilde)\colon X\to Y\) is invertible.

\section[The weighted coefficient algebra and the maps from it]{The weighted coefficient algebra\texorpdfstring{\newline}{ } and the maps from it}\label{sec:algebras}

The estimates of Section~\ref{sec:reductions} use bounded linear maps on
the coefficient space. A bounded linear map out of this space is determined
by its values on the atoms, which also give its norm bounds. The realization
theorems of Section~\ref{sec:examples} use evaluation to turn coefficient
identities into identities of functions. This section constructs the space
and specifies the weighted bound under which prescribed values on the
atoms extend to a bounded linear map, the linear lift. When these values preserve the unit and products, the
resulting map is an algebra homomorphism, the algebra lift. Evaluation
is an instance of the algebra lift, and its composition with polynomial
substitution carries polynomial identities from sequences to functions.
The categorical formulation is in Appendix~\ref{app:categorical}.
Appendix~\ref{app:lp} discusses convolution and operator bounds on
weighted \(\ell^p\) spaces.

We use the sequence space setting of~\cite[\S1.1]{lessardMJ2017},
where a numerical approximation \(\xbar\) is validated by proving the
existence of a nearby exact zero \(\xtilde\). For the examples of Section~\ref{sec:examples},
\(x(t)=\sum_{n\ge0}a_nt^n\), and matching coefficients gives \(F(a)=0\).
The finitely supported approximation \(\abar\) and the exact zero
\(\atilde\) play the roles of \(\xbar\) and \(\xtilde\) in
Theorem~\ref{thm:radii-poly}.

\subsection{The coefficient space in Lean}\label{sec:weighted-space}

\begin{definition}[Weighted sequence space, \lean{l1Weighted}]\label{def:l1Weighted}
For a fixed \(\nu\in\R, \nu>0\),
\[
  \ell^1_\nu \;:=\; \Bigl\{\, a = (a_n)_{n \ge 0} \in \R^{\N} \;\Big|\;
  \|a\|_{1,\nu} := \sum_{n=0}^{\infty} |a_n|\,\nu^n \,<\, \infty \,\Bigr\}.
\]
Each \(a\in\ell^1_\nu\) defines a power series with \(\sum_{n=0}^{\infty}|a_n|\,|z|^n\le\|a\|_{1,\nu}<\infty\) for every \(z\in\R\) with \(|z|\le\nu\)
(cf.\ \cite[\S3.1]{vdbMJLM2011}%
).
\end{definition}

The library states the space over any index set \(M\) with a positive
weight \(\omega\colon M\to(0,\infty)\), as \(\ell^1_\omega(M)\), consisting of the
family \(a=(a_m)_{m\in M}\) with \(\|a\|=\sum_m|a_m|\,\omega_m<\infty\),
and this paper specializes it to \(M=\N\) and \(\omega_n=\nu^n\). The \emph{atom}
\(\delta_n\) has coordinate \(1\) at \(n\) and \(0\) elsewhere, so
\(\|\delta_n\|=\nu^n\), and every \(a\) is the norm-convergent sum
\(a=\sum_na_n\delta_n\) of its scaled atoms.

In Lean the space is built from Mathlib's \lean{lp}. Classically
\(\ell^1\) consists of the real sequences with summable absolute values.
Mathlib states it more generally. The \(n\)-th entry of a family
lives in its own normed additive commutative group \(E(n)\), and \lean{lp}\,\(E\)\,\(1\) is the space of families \(f\) with
\(f(n)\in E(n)\) and \(\sum_n\|f(n)\|<\infty\), a Banach space when every
\(E(n)\) is complete. With \(E(n)=\R\) for every \(n\) this is the classical
\(\ell^1\). The weighted norm \(|a_n|\,\nu^n\) is obtained by letting the
\(n\)-th entry live in a copy of \(\R\) that carries the norm \(|s|\,\nu^n\)
in place of the usual absolute value:
\begin{lstlisting}
def WeightedScalar (�� : Type*) {M : Type*} (_w : M → ℝ) (_m : M) := ��

abbrev ScaledReal (ν : PosReal) := WeightedScalar ℝ (fun n : ℕ => (ν : ℝ) ^ n)
\end{lstlisting}
\lean{WeightedScalar}\:\(\Bbbk\)\:\(w\)\:\(n\) is a copy of the field
\(\Bbbk\) with the norm \(|s|\,w(n)\), and \lean{ScaledReal}\,\(\nu\) is
the copy of \(\R\) at the weight \(n\mapsto\nu^n\), so that
\lean{lp}\,(\lean{ScaledReal}\:\(\nu\))\,\(1\) is \(\ell^1_\nu\) as a Banach
space. Mathlib's \lean{lp} does not supply the Cauchy product on \(\ell^1_\nu\). The library wraps it in a
structure of its own, on which the Cauchy product of
Section~\ref{sec:monoid-algebra} is then declared:
\begin{lstlisting}
structure lpOneAlg (M : Type*) (E : M → Type*)
    [∀ m, NormedAddCommGroup (E m)] where
  toLp : lp E 1
\end{lstlisting}
An element of \lean{lpOneAlg}\:\(M\)\:\(E\) is an element of
\lean{lp}\,\(E\)\,\(1\) and nothing more. Combining this wrapper with \lean{ScaledReal} gives
the coefficient space of this paper:
\begin{lstlisting}
abbrev l1Weighted (ν : PosReal) := lpOneAlg ℕ (ScaledReal ν)
\end{lstlisting}
\lean{l1Weighted}\:\(\nu\) is \(\ell^1_\nu\) of
Definition~\ref{def:l1Weighted}. An element is a family of real numbers
\(a_n\), the \(n\)-th entry carrying the norm \(|a_n|\,\nu^n\), so that membership
in it reads \(\sum_n|a_n|\,\nu^n<\infty\) and the norm is
that sum. Completeness is inherited from Mathlib's \lean{lp} through the
wrapper:
\begin{lstlisting}
instance lpOneAlg.instCompleteSpace [∀ m, CompleteSpace (E m)] :
    CompleteSpace (lpOneAlg M E)
\end{lstlisting}
Each weighted copy of \(\R\) is complete, so this instance applies to
\(\ell^1_\nu\).

In the classical \(\ell^1_\nu\) of Definition~\ref{def:l1Weighted} the
\(n\)-th coordinate of a sequence is a real number. In the Lean space it is
an element of \(E(n)\), the copy \lean{ScaledReal}\:\(\nu\)\:\(n\) of
\(\R\). This copy is the \emph{fiber} over \(n\), and the map
\(\theta_n\colon\R\to E(n)\) sending a real number to its copy identifies
the two. The general statements of the library are proved for every space
\lean{lp}\,\(E\)\,1 whose fibers are copies of the scalars, and they need
this identification named. The class that names it is
\begin{lstlisting}
class FiberTrivialization (�� : outParam Type*) [NormedField ��]
    (M : Type*) (E : M → Type*)
    [∀ m, NormedAddCommGroup (E m)] [∀ m, NormedSpace �� (E m)] where
  equiv : ∀ m, �� ≃ₗ[��] E m
\end{lstlisting}
For a general index set \(M\), scalar field \(\Bbbk\) and fibers \(E\), the
trivialization is a linear equivalence
\(\theta_n\colon\Bbbk\simeq E(n)\) for every index. Write
\(\sigma_n=\theta_n(1)\in E(n)\), the chosen generator of the fiber. The
atom \(\delta_n\) is the sequence whose entry is \(\sigma_n\) at \(n\) and
zero otherwise. For an element \(a\) of \lean{lp}\,\(E\)\,1, write \(a(n)\in E(n)\)
for its entry and \(a_n\in\Bbbk\) for its scalar coordinate. Then
\[
 a_n=\theta_n^{-1}(a(n)),\qquad a(n)=a_n\,\sigma_n,\qquad
 \omega_n=\|\sigma_n\|=\|\delta_n\|,
\]
and the norm satisfies \(\|\theta_n(r)\|=|r|\,\omega_n\).
Section~\ref{sec:monoid-algebra} uses this trivialization to carry the
Cauchy product of scalar coordinates back to the fibers. The weight
\(\omega_n\) used in Sections~\ref{sec:free-banach}
and~\ref{sec:monoid-algebra} is \lean{FiberTrivialization.weight}. For \lean{WeightedScalar}\:\(\Bbbk\)\:\(w\)
the trivialization sends a scalar to its copy as described before, and the weight is \(w\),
as proved by \lean{WeightedScalar.}\allowbreak\lean{fiberTrivialization\_weight}. For
\lean{ScaledReal}\:\(\nu\) it is therefore \(\omega_n=\nu^n\).

\subsection{Linear maps out of the coefficient space}\label{sec:free-banach}

Throughout, \(a\in\ell^1_\nu\) has the coordinates
\(a_n=\theta_n^{-1}(a(n))\) of Section~\ref{sec:weighted-space} and is the
sum \(a=\sum_na_n\delta_n\) of its scaled atoms. Let
\(Y\) be a real Banach space. A linear map out of \(\ell^1_\nu\) into
\(Y\) is most easily specified by its values on the atoms, one vector
\(y_n\in Y\) for each index. On a finitely supported sequence this already
fixes the value, the finite sum \(\sum_na_ny_n\). Every \(a\) is the limit
of its truncations \(\sum_{n\le N}a_n\delta_n\) in the weighted norm, so
a bounded extension exists exactly when there is one constant \(C\)
such that every finitely supported \(a\) satisfies
\[
 \left\|\sum_n a_ny_n\right\|
 \le C\sum_n|a_n|\,\nu^n=C\,\|a\|_{1,\nu}.
\]
This is equivalent to the growth
condition \(\|y_n\|\le C\,\nu^n\) for every \(n\). The space
\(\ell^1_\nu\) is in this sense the completion, in the weighted norm, of
the finite linear combinations of its atoms
(Appendix~\ref{app:categorical}), and every bounded linear map out of it
is determined by its values on the atoms.

For such a family \(y\colon\N\to Y\), the extension is the bounded linear map
\(\hat y\colon\ell^1_\nu\to Y\), given by the series
\[
 \hat y(a)=\sum_{n\ge0}a_n\,y_n,\qquad \hat y(\delta_n)=y_n,\qquad
 \|\hat y(a)\|\le C\,\|a\|_{1,\nu},
\]
which converges absolutely and on a finitely supported \(a\) is the finite
sum above. Its operator norm is
\(\|\hat y\|=\sup_{a\ne0}\|\hat y(a)\|/\|a\|_{1,\nu}\), the least
constant in such a bound. We call \(\hat y\) the \emph{lift} of \(y\), and
\(y\) its \emph{column family}, by analogy with a finite matrix, whose
columns are the images of the basis vectors. Two bounded linear maps out
of \(\ell^1_\nu\) that agree on every atom are equal. This is the universal property
of \(\ell^1_\nu\) as the free Banach space on the weighted set
\((\N,\nu^\bullet)\) \cite[Thm.~3.1.1, Cor.~3.1.3]{Grilliette2012}. In the diagram,
\(\eta\) sends \(n\) to its atom \(\delta_n\), \(y\) is the
given family, and \(\hat y\) the unique bounded linear map through which
\(y\) factors:
\[
\begin{tikzcd}[column sep=large,row sep=large,ampersand replacement=\&]
 \N\arrow[r,"\eta"]\arrow[dr,"y"']
 \&\ell^1_\nu\arrow[d,"{\exists!\,\hat y}",dashed]\\
 \&Y
\end{tikzcd}
\qquad
 \|\hat y\|=\sup_{n\ge0}\frac{\|y_n\|}{\nu^n}.
\]
The dashed arrow is \lean{lpOneAlg.liftCLM}. The declarations use a general
index set \(M\), scalar field \(\Bbbk\), and fibers \(E\). Here these specialize
to \(M=\N\), \(\Bbbk=\R\), and \(E(n)=\lean{ScaledReal}\,\nu\,n\), and the following
code block's \lean{F} and \lean{v} are \(Y\) and \(y\).
\begin{lstlisting}
def liftCLM (v : M → F) (C : ℝ)
    (hv : ∀ m, ‖v m‖ ≤ C * FiberTrivialization.weight (E := E) m) :
    lpOneAlg M E →L[��] F
\end{lstlisting}
The formula is \lean{liftCLM\_apply}, and uniqueness is
\lean{continuousLinearMap\_ext}. The bound on \lean{v} is the growth
condition \(\|y_n\|\le C\,\nu^n\). Every bounded linear
map \(W\colon\ell^1_\nu\to Y\) arises this way, as the lift of its column
family \(y_n=W(\delta_n)\). For the operator estimates, it suffices to
bound the weighted column norms uniformly. The lift bound is
\lean{norm\_}\allowbreak\lean{liftCLM\_le}, and for an existing map \(W\),
\lean{norm\_le\_of\_cols} gives
\begin{equation}\label{eq:column-principle}
 \|W\delta_n\|\le C\,\nu^n\quad\text{for all }n
 \quad\implies\quad
 \|Wa\|\le C\,\|a\|\quad\text{for all }a\in\ell^1_\nu.
\end{equation}
For \(C\ge0\), this yields \(\|W\|\le C\).

For an \((N+1)\times(N+1)\) matrix \(H\), use the weighted norm
\(\|\cdot\|_{1,\nu}\) on both its domain and its codomain. The vectors
\(H\delta_j\) are its columns. Since \(\|\delta_j\|_{1,\nu}=\nu^j\), a bound
\(\|H\|\le C\) requires \(\|H\delta_j\|_{1,\nu}\le C\nu^j\) for every \(j\).
Conversely, equation~\eqref{eq:column-principle} gives \(\|H\|\le C\)
whenever these column bounds hold. The
smallest such \(C\) is therefore the largest normalized column norm,
\[
 \|H\|=\max_{0\le j\le N}\frac{\|H\delta_j\|_{1,\nu}}{\nu^j}.
\]
We also write this operator norm as \(\|H\|_{1,\nu}\).
Section~\ref{sec:reductions} uses these column bounds for its operator estimates.

\subsection{The Cauchy product and the algebra lift}\label{sec:monoid-algebra}\label{sec:space}

For evaluation of power series, the linear lift must also preserve
multiplication. On sequences, multiplication is given by the
Cauchy product.

\begin{definition}[Cauchy product, \lean{cauchyProduct}]\label{def:cauchy}
Let \(R\) be a semiring. For sequences \(f, g : \N \to R\) and
\(n\in\N\),
\[
  (f \star g)(n) \;=\; \sum_{i+j = n} f(i)\, g(j)
\]
(cf.\ \cite{kaliesKMJ2018}%
).
\end{definition}

For formal power series in \(R[[X]]\),
\[
  \sum_n (f\star g)_n X^n
  = \left(\sum_n f_n X^n\right)\left(\sum_n g_n X^n\right),
\]
so the Cauchy product is the product of formal power series, written on
the coefficients.

On \(\ell^1_\nu\) the Cauchy product of two sequences lies again in
\(\ell^1_\nu\), with
\[
 \|a\star b\|_{1,\nu}\le\|a\|_{1,\nu}\|b\|_{1,\nu},\qquad
 \delta_m\star\delta_n=\delta_{m+n},\qquad 1=\delta_0.
\]
In Lean the product is carried through the trivialization \(\theta\) of
Section~\ref{sec:weighted-space}, with \(a_m\) and \(b_n\) the scalar
coordinates and \((ab)(k)\) the entry:
\[
 (ab)(k)=\theta_k\Bigl(\sum_{m+n=k}a_m\,b_n\Bigr),\qquad 1=\delta_0.
\]
As in Section~\ref{sec:free-banach}, the code blocks of this section are
stated for a general index set \(M\), scalar field \(\Bbbk\) and fibers \(E\),
and \(\ell^1_\nu\) is their case \(M=\N\), \(\Bbbk=\R\),
\(E(n)=\lean{ScaledReal}\,\nu\,n\). The product on \(\ell_\nu^1\) is the instance
\begin{lstlisting}
instance lpOneAlg.instMul [lpOneAlgConvCompat �� M E] : Mul (lpOneAlg M E)
\end{lstlisting}
and the class it requires, \lean{lpOneAlgConvCompat}, has three fields.
For every \(k\) the sum over \(m+n=k\) converges (\lean{convSummable}),
the triple sums behind associativity converge
(\lean{tripleConvSummable}), and the norm of the \(k\)-th entry of the
product is at most \(\sum_{m+n=k}\|a(m)\|\,\|b(n)\|\), the entry bound
(\lean{norm\_conv\_le\_fiber}).
\begin{lstlisting}
class lpOneAlgConvCompat (�� : outParam Type*) [NormedField ��]
    (M : Type*) (E : M → Type*) [AddMonoid M]
    [∀ m, NormedAddCommGroup (E m)] [∀ m, NormedSpace �� (E m)]
    [FiberTrivialization �� M E] where
  convSummable : ∀ (f g : lpOneAlg M E) (k : M),
    Summable fun ab : addFiber k =>
      lpOneAlg.toRealSeq f ab.1.1 * lpOneAlg.toRealSeq g ab.1.2
  tripleConvSummable : ∀ (f g h : lpOneAlg M E) (x : M),
    AddTripleConvolutionSummable
      (lpOneAlg.toRealSeq f) (lpOneAlg.toRealSeq g) (lpOneAlg.toRealSeq h) x
  norm_conv_le_fiber : ∀ (f g : lpOneAlg M E) (k : M),
    ‖FiberTrivialization.equiv (E := E) k
      (addRingConvolution (lpOneAlg.toRealSeq f)
        (lpOneAlg.toRealSeq g) k)‖ ≤
    ∑' ab : addFiber k, ‖f ab.1.1‖ * ‖g ab.1.2‖
\end{lstlisting}
The entry bound writes the convolution with Mathlib's
\lean{addRingConvolution}, which on \(\N\) equals \lean{cauchyProduct}
(\lean{cauchyProduct.eq\_addRingConvolution}). Summed over \(k\), it is
the bound \(\|a\star b\|_{1,\nu}\le\|a\|_{1,\nu}\|b\|_{1,\nu}\) displayed
above. On \(\N\) the pair and triple sums are finite, and the triangle
inequality with \(\nu^{m+n}=\nu^m\nu^n\) gives the entry bound, so
\(\ell^1_\nu\) satisfies \lean{lpOneAlgConvCompat} and carries the product.
Mathlib's \lean{lp} supplies the normed
space structure, and with the product, \(\ell^1_\nu\) is a unital real
Banach algebra:
\begin{lstlisting}
instance lpOneAlg.instNormedCommRing : NormedCommRing (lpOneAlg M E)

instance lpOneAlg.instNormedAlgebra : NormedAlgebra �� (lpOneAlg M E)
\end{lstlisting}

\paragraph{The algebra lift.}
Let \(Y\) now be a unital real Banach algebra and \(y\colon\N\to Y\) a
family with the growth condition \(\|y_n\|\le C\,\nu^n\). The atoms
multiply by \(\delta_m\star\delta_n=\delta_{m+n}\), with unit
\(\delta_0\). Compared with the linear lift of
Section~\ref{sec:free-banach}, the additional conditions ask \(y\) to carry
this unit and these products to the unit and the products of \(Y\):
\[
 y_0=1_Y,\qquad y_{m+n}=y_m\,y_n\quad(m,n\in\N).
\]
The linear lift then fits the same triangle, now as an algebra homomorphism:
\[
\begin{tikzcd}[column sep=large,row sep=large,ampersand replacement=\&]
 (\N,+)\arrow[r,"\eta"]\arrow[dr,"y"']
 \&\ell^1_\nu\arrow[d,"{\exists!\,\hat y}",dashed]\\
 \&(Y,\cdot)
\end{tikzcd}
\]
The dashed arrow is \lean{lpOneAlg.liftAlgHom}:
\begin{lstlisting}
def liftAlgHom (b : M → B) (C : ℝ)
    (hb : ∀ m, ‖b m‖ ≤ C * FiberTrivialization.weight (E := E) m)
    (hb0 : b 0 = 1) (hbmul : ∀ m n, b (m + n) = b m * b n) :
    lpOneAlg M E →A[��] B
\end{lstlisting}
Its underlying linear map is the lift of Section~\ref{sec:free-banach},
and it sends \(\delta_n\) to \(y_n\). Here \lean{B} and \lean{b} are
\(Y\) and \(y\), and the conditions \lean{hb0} and \lean{hbmul} are the
two displayed conditions. For each fixed factor, both sides of
\(\hat y(a\star b)=\hat y(a)\,\hat y(b)\) are bounded linear in the other
factor. They agree on pairs of atoms by \(y_{m+n}=y_my_n\), so the
uniqueness theorem of Section~\ref{sec:free-banach}, applied in each
factor, gives the identity (\lean{liftCLM\_mul}). Uniqueness of the
algebra lift is \lean{continuousAlgHom\_ext}. Every continuous algebra
homomorphism \(\psi\colon\ell^1_\nu\to Y\) arises this way, as the lift
of the family \(y_n=\psi(\delta_n)\), which satisfies the two conditions
because \(\psi\) is unital and multiplicative. Appendix~\ref{app:categorical}
describes this lift through the universal property of the free Banach
algebra on a weighted monoid.

\subsection{Evaluation}\label{sec:evaluations}

In this section we identify a polynomial \(\sum_{n=0}^N a_nx^n\in\R[x]\) with the finite
sequence \(\sum_{n=0}^N a_n\delta_n\). Addition and the Cauchy product
then agree with polynomial addition and multiplication, with \(x\)
corresponding to \(\delta_1\). For a sequence \(a=(a_n)\in\ell^1_\nu\), the associated
power series is \(\sum_{n\ge0}a_nx^n\). Its partial sums correspond to
the truncations of \(a\), which converge to \(a\) in the weighted norm.

Let \(Y\) be a unital real Banach algebra. Polynomial evaluation at
\(x\in Y\) sends \(\delta_n\) to \(x^n\). If \(\|x^n\|\le C\,\nu^n\) for
all \(n\), the algebra lift of Section~\ref{sec:monoid-algebra}, applied
to the family \(n\mapsto x^n\), extends polynomial evaluation to the
continuous algebra homomorphism \(\operatorname{ev}_x\colon\ell^1_\nu\to Y\),
given by \(\operatorname{ev}_x(a)=\sum_na_nx^n\). In the diagram, \(\delta_1\)
is the generator, \(x\) the given point in \(Y\) with \(\|x^n\|\le C\,\nu^n\), and
\(\operatorname{ev}_x\) the unique continuous unital algebra homomorphism
through which \(x\) factors:
\[
\begin{tikzcd}[column sep=large,row sep=large,ampersand replacement=\&]
 \{1\}\arrow[r,"\delta_1"]\arrow[dr,"x"']
 \&\ell^1_\nu\arrow[d,"{\exists!\,\operatorname{ev}_x}",dashed]\\
 \&Y
\end{tikzcd}
\]
The dashed arrow is \lean{lpOneAlg.geomAeval}:
\begin{lstlisting}
def geomAeval (x : B) (C : ℝ) (hx : ∀ n, ‖x ^ n‖ ≤ C * (ν : ℝ) ^ n) :
    lpOneAlg ℕ (geomFiber �� ν) →A[��] B :=
  liftAlgHom (fun n => x ^ n) C (pow_le_geomFiber_weight �� ν x C hx)
    (pow_zero x) (pow_add x)
\end{lstlisting}
This is \lean{liftAlgHom} at the family \(y_n=x^n\). The growth condition
\lean{hb} expresses \lean{hx}, since the weight is \(\nu^n\), and \lean{hb0} and
\lean{hbmul} are \lean{pow\_zero} and \lean{pow\_add}. Its fiber family is defined by
\begin{lstlisting}
abbrev geomFiber (�� : Type*) [NormedField ��] (ν : PosReal) :=
  WeightedScalar �� (fun n : ℕ => (ν : ℝ) ^ n)
\end{lstlisting}
For the real coefficients used in this paper, \(\Bbbk=\R\), 
so \lean{geomFiber}\:\(\R\)\:\(\nu\) is \lean{ScaledReal}\:\(\nu\) of
Section~\ref{sec:weighted-space}, and the domain is
\lean{l1Weighted}\,\(\nu\). The resulting continuous unital
algebra homomorphism sends \(\delta_1\) to \(x\), and
\lean{geomAlgHom\_ext} gives its uniqueness. Conversely, every continuous
unital algebra homomorphism \(\psi\colon\ell^1_\nu\to Y\) is evaluation at
\(x=\psi(\delta_1)\), since \(\psi(\delta_n)=x^n\). Under the growth condition,
\(\|\operatorname{ev}_x\|=\sup_n\|x^n\|/\nu^n\le C\).

Every \(x\) with \(\|x\|\le\nu\) satisfies the growth condition with
\(C=\max\{1,\|1\|\}\). Evaluation at a real number is the case \(Y=\R\), \(C=1\).

\begin{theorem}[Cauchy-product evaluation, \lean{eval\_mul}]\label{thm:mertens}
Fix \(x\in[-\nu,\nu]\). The evaluation map
\[
  \lean{eval}(\cdot,x)\colon \ell^1_\nu\to\R,\qquad
  a\longmapsto\sum_{n\ge0}a_nx^n,
\]
is well defined, since the series converges absolutely, and it is multiplicative:
\[
  \lean{eval}(a\star b,x)
  = \lean{eval}(a,x)\,\lean{eval}(b,x)
  \qquad(a,b\in\ell^1_\nu).
\]
\end{theorem}

For \(x\in\R\) with \(|x|\le\nu\), the lift \(\lean{aeval\_of\_norm\_le}\,\nu\,x\)
sends \(\delta_1\) to \(x\), and its underlying function is
\(\lean{eval}(\cdot,x)\), with norm at most \(1\), that is
\(|\lean{eval}(a,x)|\le\|a\|_{1,\nu}\).

\subsection{Free commutative algebra and polynomials}\label{sec:free}

For the Taylor instances of \textbf{P2}, the polynomial components are evaluated
on real vectors for the differential equation and on coefficient sequences
for the map \(F\). The certificate computes coefficients of the same
expressions from rational arrays. The library represents each component
once and uses Mathlib's polynomial ring \MvPoly{} to relate these
evaluations.

\begin{definition}[\(\MvPoly\)]\label{def:mvpoly}
For a commutative ring \(R\) and index set \(\sigma\), the ring
\(\MvPoly\,\sigma\,R\) is the polynomial ring
\(R[X_i : i \in \sigma]\) of finite formal \(R\)-linear combinations
of monomials \(\prod_{i \in \sigma} X_i^{n_i}\) with finitely many
\(n_i\) nonzero.
\end{definition}

\begin{definition}[Universal property, \lean{aeval}]\label{def:aeval}
Given a commutative ring \(R\), an index set \(\sigma\), a commutative
\(R\)-algebra \(S\), and a function \(a\colon \sigma \to S\), there is a
unique \(R\)-algebra homomorphism
\(\operatorname{ev}_a\colon\MvPoly\,\sigma\,R\to S\) through which \(a\)
factors:
\[
\begin{tikzcd}[column sep=large,row sep=large,ampersand replacement=\&]
 \sigma\arrow[r,"{i\mapsto X_i}"]\arrow[dr,"a"']
 \&R[X_i : i\in\sigma]\arrow[d,"{\exists!\,\operatorname{ev}_a}",dashed]\\
 \&S
\end{tikzcd}
\]
\end{definition}
The dashed arrow is Mathlib's \lean{MvPolynomial.aeval}~\(a\):
\begin{lstlisting}
def aeval (f : σ → S₁) : MvPolynomial σ R →ₐ[R] S₁
\end{lstlisting}
with \(\lean{f}=a\) and \(\lean{S}_1=S\). The generator equation is \(\operatorname{ev}_a\,(X_i)=a(i)\). For every \(R\)-algebra homomorphism
\(h\colon S\to S'\),
\begin{equation}\label{eq:aeval-natural}
  h\circ\operatorname{ev}_a \;=\; \operatorname{ev}_{h\circ a},
\end{equation}
so polynomial substitution commutes with every algebra homomorphism out of
the codomain, as stated by Mathlib's \lean{MvPolynomial.comp\_aeval}.

\paragraph{Differentiation.}
Take \(R=\R\) and \(S=\ell^1_\nu\) in Definition~\ref{def:aeval}.
Fix a polynomial \(P\in\R[X]\) in one variable. Varying
\(a\in\ell^1_\nu\) gives the polynomial map
\[
 \varphi\colon\ell^1_\nu\to\ell^1_\nu,\qquad \varphi(a)=\operatorname{ev}_a(P).
\]

The norm bounds in Theorem~\ref{thm:radii-poly} require Fr\'echet
derivatives. For the polynomial map \(\varphi\), this derivative is
computed from the formal derivative \(\partial P\).

For \(P=X^2\), \(\varphi(a)=a\star a\). The map
\(b\mapsto 2(a\star b)\) is bounded and linear. For nonzero
\(b\in\ell^1_\nu\), the remainder satisfies
\[
 \frac{\|\varphi(a+b)-\varphi(a)-2(a\star b)\|}{\|b\|}
 =\frac{\|b\star b\|}{\|b\|}
 \le\|b\|\longrightarrow0\qquad(b\to0).
\]
Thus the Fr\'echet derivative of \(\varphi\) at \(a\) is the left action of
the element \(2a=\operatorname{ev}_a(\partial P)\), since \(\partial P=2X\).

For any polynomial \(P\) in one variable, Mathlib's
\lean{MvPolynomial.pderiv} forms the formal derivative \(\partial P\). In every normed commutative \(\R\)-algebra
\(Y\), the library's \lean{MvPolynomial.}\allowbreak\lean{hasFDerivAt\_aeval}
identifies the derivative
of \(a\mapsto\operatorname{ev}_a(P)\) at \(a\in Y\) with the left action of the element
\(\operatorname{ev}_a(\partial P)\), the linear map
\[
 b\longmapsto\operatorname{ev}_a(\partial P)\,b\qquad(b\in Y).
\]
For \(Y=\ell^1_\nu\), the library writes this left action as
\lean{l1Weighted.leftMul} applied to \(\operatorname{ev}_a(\partial P)\). For \(L\) variables, take
\(\sigma=\{0,\dots,L-1\}\) and let
\(a=(a_0,\dots,a_{L-1})\in(\ell^1_\nu)^L\).
The library states the derivative formula for polynomials with rational
coefficients as \lean{fderiv\_aeval}:
\begin{lstlisting}
theorem fderiv_aeval (p : MvPolynomial (Fin L) ℚ) (a : Fin L → l1Weighted ν) :
    fderiv ℝ (fun x => MvPolynomial.aeval x p) a =
      ∑ i : Fin L,
        (l1Weighted.leftMul (MvPolynomial.aeval a (MvPolynomial.pderiv i p))).comp
          (proj (R := ℝ) (φ := fun _ : Fin L => l1Weighted ν) i)
\end{lstlisting}

\paragraph{Substitution and evaluation.}
Substituting coefficient sequences in a polynomial and evaluating the
resulting series gives the identity of Section~\ref{sec:intro}.
Write \(\R[X]=\R[X_0,\dots,X_{L-1}]\). For \(a\in(\ell^1_\nu)^L\), let
\(x\) be a point of a commutative complete normed \(\R\)-algebra \(Y\) with
\(\|x^n\|\le C\,\nu^n\) for every \(n\).
Definition~\ref{def:aeval} induces the horizontal and lower diagonal
arrows from \(a\) and \(\operatorname{ev}_x\circ a\), respectively. The
right vertical arrow is evaluation at \(x\) from
Section~\ref{sec:evaluations}.
\[
\begin{tikzcd}[column sep=large,row sep=large,ampersand replacement=\&]
 \{0,\dots,L-1\}\arrow[d,"{i\mapsto X_i}"']\arrow[dr,"a"]
 \&\\
 \R[X]\arrow[r,"\operatorname{ev}_a"]\arrow[dr,"\operatorname{ev}_{\operatorname{ev}_x\circ a}"']
 \&\ell^1_\nu\arrow[d,"\operatorname{ev}_x"]\\
 \&Y
\end{tikzcd}
\]
For every \(P\in\R[X]\), the lower triangle reads
\[
 \operatorname{ev}_x\bigl(\operatorname{ev}_a(P)\bigr)
 =\operatorname{ev}_{\operatorname{ev}_x\circ a}(P),
\]
which is \eqref{eq:aeval-natural} at \(h=\operatorname{ev}_x\). Since
\(\operatorname{ev}_a\) is a homomorphism, \(\operatorname{ev}_a(P)\) is \(P\) with
\(a_i\) substituted for \(X_i\), written \(P(a_0,\dots,a_{L-1})\). In this
notation the display is the identity of Section~\ref{sec:intro}. On the left, \(P\) is evaluated on
the sequences and the resulting series is summed at \(x\). On
the right, each series \(a_i\) is summed at \(x\) and \(P\) is evaluated at the
values \(\operatorname{ev}_x(a_0),\dots,\operatorname{ev}_x(a_{L-1})\). For
\(Y=\R\), these values form a vector in \(\R^L\), and the right side is
the polynomial evaluation used in the differential equation. The construction of \(\operatorname{ev}_a\) uses no norm or completeness,
since the \(X_i\) generate algebraically. Analysis enters only through
the existence of \(\operatorname{ev}_x\)
(Section~\ref{sec:evaluations}). Appendix~\ref{app:categorical} reads
the diagram with its categories and the forgetful functor written.

\par\addvspace{3.25ex plus 1ex minus .2ex}
\noindent\begin{minipage}{\linewidth}
\paragraph{Computable representation.}
Mathlib implements \MvPoly{} noncomputably, so a certificate cannot
evaluate its terms by computation. The library therefore stores a polynomial
as a term of the inductive type \lean{CompPoly}~\(L\), a computable
polynomial expression over \(\Q\) in the variables \(X_0,\ldots,X_{L-1}\):
\begin{lstlisting}
inductive CompPoly (L : ℕ) where
  | C : ℚ → CompPoly L
  | X : Fin L → CompPoly L
  | add : CompPoly L → CompPoly L → CompPoly L
  | sub : CompPoly L → CompPoly L → CompPoly L
  | mul : CompPoly L → CompPoly L → CompPoly L
  | neg : CompPoly L → CompPoly L
  | smul : ℚ → CompPoly L → CompPoly L
\end{lstlisting}
\end{minipage}\par
\Needspace{11\baselineskip}
We call a term \(p\) of type \lean{CompPoly}~\(L\) an \emph{expression}. Its
\emph{translation} is the polynomial \(P=p.\lean{toMvPoly}\) that it represents
in \(\Q[X]=\Q[X_0,\dots,X_{L-1}]\).
In the certificates of Section~\ref{sec:examples}, the numerical
approximation \(\abar\) is stored as finite arrays of rational numbers.
For each variable \(X_i\), the evaluator receives one array \(c_i\) of
type \lean{Array}~\(\Q\). The evaluator \lean{evalCoeff} below reads \(c_i\) as a
sequence, taking value \(0\) beyond the end of the array.
Casting this sequence to \(\R\) gives a finitely supported
\(a_i\in\ell^1_\nu\). In this diagram \(c=(c_i)\) denotes the input
arrays and \(a=(a_i)\) their real coefficient sequences.
In the polynomial-zero family \(a_0=\abar\), and in the Taylor family
\(a=\abar\). The certificate computes with \(c\) to obtain
coefficients of the expression evaluated at \(a\). The diagram below connects this computation with evaluation at \(a\). Its
identities hold for every expression, so a certificate for a new polynomial
needs no new proof of them.
\[
\begin{tikzcd}[column sep=huge,row sep=huge,ampersand replacement=\&]
 \&\text{expressions}
   \arrow[dl,"{\substack{\text{coefficients at }c\\\text{(computable)}\\\text{\lean{evalCoeff}}}}"']
   \arrow[d,"{\substack{\text{translation}\\\text{\lean{toMvPoly}}}}" pos=0.70]
   \arrow[dr,"{\substack{\text{evaluation at }a\\\text{(noncomputable)}\\\text{\lean{evalBanach}}}}"]\&\\
 (\N\to\Q)\arrow[dr,"\text{real cast}"']
 \&\Q[X]\arrow[l,"{\substack{P\mapsto\operatorname{coeff}(P(c))\\\text{\lean{mvPolyCoeff}}}}"']
          \arrow[r,"{\substack{P\mapsto P(a)\\\text{\lean{aeval}}}}"]
 \&\ell^1_\nu\arrow[dl,"{\substack{\text{coefficients}\\\text{\lean{toSeq}}}}"]\\
 \&(\N\to\R)\&
\end{tikzcd}
\]
The two outer arrows evaluate the expression independently, by recursion on
its syntax, in rational sequences on the left and in \(\ell^1_\nu\) on the
right. The middle arrow translates it to \(P\). The middle row then
evaluates the same \(P\) in two algebras by Definition~\ref{def:aeval},
\[
 \Q[[t]]\ \xleftarrow{\ \operatorname{ev}_c\ }\ \Q[X]\ \xrightarrow{\ \operatorname{ev}_a\ }\ \ell^1_\nu,
 \qquad \operatorname{ev}_c\colon X_i\mapsto\sum_n c_i(n)\,t^n,\qquad
 \operatorname{ev}_a\colon X_i\mapsto a_i,
\]
where \(\Q[[t]]=\lean{PowerSeries}\,\Q\). The diagram's left arrow out of
\(\Q[X]\) is \(\operatorname{ev}_c\) followed by \(\operatorname{coeff}\) (\lean{PowerSeries.coeff}),
which reads a series off as its coefficient sequence in \((\N\to\Q)\). The
direct evaluators agree with
\(\operatorname{coeff}\circ\operatorname{ev}_c\circ\lean{toMvPoly}\) on the
left and \(\operatorname{ev}_a\circ\lean{toMvPoly}\) on the right. The coefficientwise real
cast preserves sums and products and sends the coefficients \(c_i(n)\)
to those of \(a_i\), giving the bottom comparison.

An example first shows what each arrow does, and the Lean declarations that
construct the arrows follow. Take \(L=1\), the expression \(p=X_0\cdot X_0\) of Section~\ref{sec:sqrt},
and the array \(c_0=(2,1)\), read as the sequence \((2,1,0,\dots)\) with
real cast \(a_0=(2,1,0,\dots)\in\ell^1_\nu\). For these inputs, the seven
arrows give:
\[
\begin{array}{@{}lll@{\;}l@{}}
 \text{translation} & p & \longmapsto P=X_0^2 & \in\Q[X_0]\\[4pt]
 \text{coefficients at }c & p & \longmapsto (4,4,1,0,\dots) & \in(\N\to\Q)\\
 P\mapsto\operatorname{coeff}(P(c)) & X_0^2 & \longmapsto \operatorname{coeff}\bigl((2+t)^2\bigr)=(4,4,1,0,\dots) & \in(\N\to\Q)\\[4pt]
 \text{evaluation at }a & p & \longmapsto a_0\star a_0 & \in\ell^1_\nu\\
 P\mapsto P(a) & X_0^2 & \longmapsto a_0\star a_0 & \in\ell^1_\nu\\[4pt]
 \text{real cast} & (4,4,1,0,\dots) & \longmapsto (4,4,1,0,\dots) & \in(\N\to\R)\\
 \text{coefficients} & a_0\star a_0 & \longmapsto (4,4,1,0,\dots) & \in(\N\to\R)
\end{array}
\]
The coefficient at index \(n=1\) is the Cauchy sum \(2\cdot1+1\cdot2=4\) on the
left and in \(\ell^1_\nu\), and \((2+t)^2=4+4t+t^2\) in the middle. The paths
meet in four equalities, each an application of a theorem. On the left
triangle, \(p.\lean{evalCoeff}\,c=\operatorname{coeff}(P(c))\), both
\((4,4,1,0,\dots)\) (\lean{evalCoeffSeq\_eq\_}\allowbreak\lean{mvPolyCoeff}).
On the right triangle, \(p.\lean{evalBanach}\,a=P(a)\), both
\(a_0\star a_0\) (\lean{compPoly\_evalBanach\_eq\_aeval}). On the bottom
square, the real cast of \(\operatorname{coeff}(P(c))\) is
\(\lean{toSeq}(P(a))\), both \((4,4,1,0,\dots)\), under the matching input
\(\lean{toSeq}(a_0)=(2,1,0,\dots)\) (\lean{toSeq\_aeval\_eq\_mvPolyCoeff}).
Around the outer square, the real cast of \(p.\lean{evalCoeff}\,c\) is
\(\lean{toSeq}(p.\lean{evalBanach}\,a)\), both \((4,4,1,0,\dots)\)
(\lean{toSeq\_evalBanach\_eq\_evalCoeff}). The structures that \lean{aeval}
needs, a commutative \(\Q\)-algebra, on \(\lean{PowerSeries}\,\Q\) and on
\(\ell^1_\nu\), are declared once as instances and found by Lean
automatically. For
\(\ell^1_\nu\) the instance is the algebra of Section~\ref{sec:monoid-algebra} with scalars
restricted from \(\R\) to \(\Q\).

The translation \lean{toMvPoly} has one clause per constructor:
\begin{lstlisting}
noncomputable def toMvPoly : CompPoly L → MvPolynomial (Fin L) ℚ
  | .C r => MvPolynomial.C r
  | .X i => MvPolynomial.X i
  | .add p q => p.toMvPoly + q.toMvPoly
  | .sub p q => p.toMvPoly - q.toMvPoly
  | .mul p q => p.toMvPoly * q.toMvPoly
  | .neg p => -p.toMvPoly
  | .smul r p => MvPolynomial.C r * p.toMvPoly
\end{lstlisting}
Different expressions, such as \(X_0\cdot X_0\) and \(X_0\cdot X_0+0\),
can have the same translation. Equality of translations is coefficientwise
(\lean{MvPolynomial.ext}). By the universal property of \(\Q[X]\)
(Definition~\ref{def:aeval}), this is equality of the expressions up to the
laws of a commutative \(\Q\)-algebra.

\Needspace{16\baselineskip}
The outer left arrow is computable because \lean{evalCoeffSeq} recurses on the expression
and uses only finitely many rational operations for the coefficient at each index \(n\).
At a product it sums the \(n+1\) terms indexed by \(i+j=n\).\footnote{Over
\(\N\), the input sequences need not have finite support to have finitely many terms at index \(n\) of their Cauchy product. Over \(\Z\), the
fiber \(i+j=n\) is infinite. The library's
\lean{evalLaurentCoeffSeq} computes finite sums using supplied support
bounds, and its correctness theorem assumes that the inputs vanish outside
those bounds.}
The array version \lean{evalCoeff} reads the arrays as sequences:
\begin{lstlisting}
def evalCoeffSeq : CompPoly L → (Fin L → ℕ → ℚ) → ℕ → ℚ
  | .C r, _, n => if n = 0 then r else 0
  | .X i, seqs, n => seqs i n
  | .add p q, seqs, n => evalCoeffSeq p seqs n + evalCoeffSeq q seqs n
  | .sub p q, seqs, n => evalCoeffSeq p seqs n - evalCoeffSeq q seqs n
  | .mul p q, seqs, n =>
      cauchyProduct (evalCoeffSeq p seqs) (evalCoeffSeq q seqs) n
  | .neg p, seqs, n => -(evalCoeffSeq p seqs n)
  | .smul r p, seqs, n => r * evalCoeffSeq p seqs n

def evalCoeff (p : CompPoly L) (arrs : Fin L → Array ℚ) : ℕ → ℚ :=
  p.evalCoeffSeq fun i n => (arrs i).getD n 0
\end{lstlisting}

\Needspace{14\baselineskip}
The outer right arrow \lean{evalBanach} substitutes \(a_i\) for \(X_i\)
directly in the expression. Its multiplication here is the Cauchy product. The
definition works in any commutative \(\R\)-algebra, written \lean{R} below:
\begin{lstlisting}
noncomputable def evalBanach {R : Type*} [CommRing R] [Algebra ℝ R]
    (p : CompPoly L) (a : Fin L → R) : R :=
  match p with
  | .C r => algebraMap ℝ R ((r : ℚ) : ℝ)
  | .X i => a i
  | .add p q => p.evalBanach a + q.evalBanach a
  | .sub p q => p.evalBanach a - q.evalBanach a
  | .mul p q => p.evalBanach a * q.evalBanach a
  | .neg p => -(p.evalBanach a)
  | .smul r p => ((r : ℚ) : ℝ) • p.evalBanach a
\end{lstlisting}

Sections~\ref{sec:poly-derivatives} and~\ref{sec:taylor-services} use
the array form of the coefficient identity,
\lean{toSeq\_}\allowbreak\lean{evalBanach\_}\allowbreak\lean{eq\_evalCoeff}. The polynomial-zero certificate of
Section~\ref{sec:sqrt} uses its sequence form for the stored coefficients.
For every expression \(p\), if each \(\lean{toSeq}(a_i)\) is the real cast
of \(c_i\) read as a sequence, then
\[
 \bigl(p.\lean{evalCoeff}\,c\,n:\R\bigr)
 = \bigl(p.\lean{evalBanach}\,a\bigr)_n .
\]

Factoring through \(\Q[X]\), rather than proving the coefficient identity
directly by induction on expressions, lets the proofs for \lean{evalBanach}
reuse Mathlib's substitution laws and formal partial derivatives. The
comparison \lean{evalBanach\_eq\_aeval\_map} identifies
expression evaluation in any commutative \(\R\)-algebra with evaluation of
\(P\) after casting its coefficients from \(\Q\) to \(\R\).

Applying \eqref{eq:aeval-natural} with \(h=\operatorname{ev}_t\), the algebra
lift of Section~\ref{sec:evaluations}, gives the substitution identity for
expressions (\lean{eval\_evalBanach}) for \(|t|\le\nu\),
\[
 \operatorname{ev}_t\bigl(p.\lean{evalBanach}\,a\bigr)
 =p.\lean{evalBanach}\,(\operatorname{ev}_t\circ a).
\]
Section~\ref{sec:taylor-ivp} uses this identity with termwise differentiation
to verify the differential equation for \(|t|<\nu\).

Through the same comparison, the polynomial derivative formula gives the
Fr\'echet derivative of \(a\mapsto p.\lean{evalBanach}\,a\) required by the
norm bounds of Theorem~\ref{thm:radii-poly}
(\lean{hasFDerivAt\_}\allowbreak\lean{evalBanach}). For each \(i\),
\lean{CompPoly.pderiv} computes an expression \(\partial_i p\) whose
translation is \(\partial_i P\) (\lean{pderiv\_toMvPoly}). In the \(i\)-th input
variable, the partial derivative at \(a\) is the left action of the element
\((\partial_i p).\lean{evalBanach}\,a\). For
\(p=X_0\cdot X_0\), \(\partial_0 p=1\cdot X_0+X_0\cdot1\), which
translates to \(2X_0\). At the same input the coefficient at index \(n=1\) is
now \(2\), the coefficient of \(2a_0=1\star a_0+a_0\star 1\).

\section{Operator estimates}\label{sec:reductions}

Theorem~\ref{thm:radii-poly} asks for four norm bounds, while a
certificate records finitely many rational numbers. The library states each
estimate once, as a lemma that reduces a norm bound to computations on
finitely many coefficients and on polynomial expressions.

\subsection{Coefficient maps}\label{sec:fun-to-seq}

Each family certifies a zero of a coefficient map. For the polynomial zero,
Theorem~\ref{thm:radii-poly} is applied to the coefficient map itself. For
the initial value problems the approximate derivative \(A^\dagger\) is unbounded on the
coefficient space. The map \(G(a)=A(F(a))\) applies \(A\) to the coefficients
of \(F(a)\), and the theorem is applied to
\(G\colon(\ell^1_\nu)^L\to(\ell^1_\nu)^L\), with the identity in the role of the theorem's \(A\)
and the bounded operator \(AA^\dagger\) of Section~\ref{sec:finite-tail} in
the role of its \(A^\dagger\).

\paragraph{Polynomial zero.}
For \(x^2-\lambda=0\) at \(\lambda_0\), put
\(t=\lambda-\lambda_0\) and \(x(t)=\sum_{n\ge0}a_nt^n\). The coefficient map
is
\[
 F(a)=a\star a-c,\qquad c=(\lambda_0,1,0,0,\ldots),
\]
and the derivative \(DF(a)\) is the left action of the element \(2a\). Its
action is \(DF(a)h=2(a\star h)\), and in coordinates
\[
 DF(a)_{nk}=\begin{cases}2a_{n-k},&k\le n,\\0,&k>n,\end{cases}
 \qquad
 DF(a)=2\begin{pmatrix}
 a_0 & 0 & 0 & \cdots\\
 a_1 & a_0 & 0 & \cdots\\
 a_2 & a_1 & a_0 & \cdots\\
 \vdots & \vdots & \vdots & \ddots
 \end{pmatrix}.
\]

\paragraph{Taylor initial value problem.}
For initial value problems, substituting the Taylor series in the differential
equation gives a recurrence involving the derivative shift \((n+1)a_{n+1}\)
\cite[\S4.1, (35)]{kaliesKMJ2018}. Let \(L\ge1\), let \(f\colon\R^L\to\R^L\) be a polynomial vector field with rational coefficients, and
let \(x_0\in\R^L\). For \(x(t)_l=\sum_{n\ge0}(a_l)_nt^n\), with \(a_{l,n}=(a_l)_n\), let
\(\varphi_l(a)\) be the component \(f_l\) evaluated at \(a\in(\ell^1_\nu)^L\)
as in Definition~\ref{def:aeval}, with \(f_l\) stored as the expression
\(q_l\) of \lean{CompPoly}~\(L\). The
coefficient map of \(\dot x=f(x)\), \(x(0)=x_0\), is
\begin{equation}\label{eq:ivp-coeff}
 \begin{aligned}
 F(a)_{l,0}&=a_{l,0}-(x_0)_l,\\
 F(a)_{l,n+1}&=(n+1)a_{l,n+1}-\varphi_l(a)_n.
 \end{aligned}
\end{equation}
Because of the factor \(n+1\), \(a\in(\ell^1_\nu)^L\) does not imply
\(F(a)\in(\ell^1_\nu)^L\). The library returns \(F(a)\) as a raw coefficient
family (\lean{ivpCoeffs}, a function \(\lean{Fin}\,L\to\N\to\R\)),
and only the preconditioned map \(G\) lands in the Banach space.\footnote{For
the realization, a second weighted space, with weights
\(\omega_n=\nu^{n+1}/(n+1)\), contains the termwise differentiated series,
and \(\frac{d}{dt}\lean{eval}(a,t)=\sum_{n\ge0}(n+1)a_{n+1}t^n\) for
\(|t|<\nu\) (Section~\ref{sec:taylor-ivp}).}

\subsection{Polynomial calculus and expression constants}\label{sec:poly-derivatives}

A formal derivative is useful only once it is identified with the
Fr\'echet derivative that Theorem~\ref{thm:radii-poly} sees. For a
\lean{CompPoly} expression this identification is
\lean{hasFDerivAt\_evalBanach} of Section~\ref{sec:free}. In the \(j\)-th
input variable, the partial derivative of \(a\mapsto q_l.\lean{evalBanach}\,a\) at \(\abar\) is
the left action of the element \((\partial_jq_l).\lean{evalBanach}\,\abar\). Let
\(c=(c_i)_{0\le i<L}\) be the rational arrays whose real casts give the
coefficients of \(\abar\in(\ell^1_\nu)^L\). The coefficients of that element are
\[
 \bigl((\partial_jq_l).\lean{evalBanach}\,\abar\bigr)_k=\bigl((\partial_jq_l).\lean{evalCoeff}\,c\,k:\R\bigr)
 \qquad(k\in\N).
\]

Constants read off an expression bound norms. Let \(\mathcal A\) be a normed
commutative ring with a normed \(\R\)-algebra structure and
\(\|1_{\mathcal A}\|=1\). For \(q\colon\lean{CompPoly}\,L\),
\(a,b\colon\lean{Fin}\,L\to\mathcal A\) and radii
\(R\colon\lean{Fin}\,L\to\R\) with \(\|a(i)\|,\|b(i)\|\le R_i\) for every
\(i\), with the supremum norm on \(\lean{Fin}\,L\to\mathcal A\),
\begin{align*}
 \|q.\lean{evalBanach}\,a\|&\le q.\lean{normBound}\,R,\\
 \|q.\lean{evalBanach}\,a-q.\lean{evalBanach}\,b\|
 &\le q.\lean{lipschitzBound}\,R\,\|a-b\|
\end{align*}
Both constants are computed by recursion on the expression with triangle
inequalities. For subexpressions \(p\) and \(q\) and a rational constant \(r\),
with \(\operatorname{nb}(p)=p.\lean{normBound}\,R\) and
\(\operatorname{lb}(p)=p.\lean{lipschitzBound}\,R\),
\[
\begin{aligned}
 \operatorname{nb}(r)&=|r|,\quad \operatorname{nb}(X_i)=R_i, &\qquad
 \operatorname{lb}(r)&=0,\quad \operatorname{lb}(X_i)=1,\\
 \operatorname{nb}(p\pm q)&=\operatorname{nb}(p)+\operatorname{nb}(q), &
 \operatorname{lb}(p\pm q)&=\operatorname{lb}(p)+\operatorname{lb}(q),\\
 \operatorname{nb}(p\cdot q)&=\operatorname{nb}(p)\operatorname{nb}(q), &
 \operatorname{lb}(p\cdot q)&=\operatorname{lb}(p)\operatorname{nb}(q)+\operatorname{nb}(p)\operatorname{lb}(q),
\end{aligned}
\]
and both are unchanged by \(-p\) and multiplied by \(|r|\) under
\(r\cdot p\). They see no cancellation, so two expressions of the same
polynomial can have different constants. The derivative constants
\[
\begin{aligned}
 q.\lean{derivativeBound}\,R&=\sum_{i=0}^{L-1}(\partial_iq).\lean{normBound}\,R,\\
 q.\lean{derivativeLipschitzBound}\,R&=\sum_{i=0}^{L-1}(\partial_iq).\lean{lipschitzBound}\,R
\end{aligned}
\]
bound the derivative of \(\Phi_q(a)=q.\lean{evalBanach}\,a\) on \((\ell^1_\nu)^L\)
and its variation. For \(\|a_i\|,\|b_i\|\le R_i\) and every direction \(h\),
\[
\begin{aligned}
 \|D\Phi_q(a)\,h\|&\le q.\lean{derivativeBound}\,R\,\|h\|,\\
 \bigl\|\bigl(D\Phi_q(b)-D\Phi_q(a)\bigr)h\bigr\|
 &\le q.\lean{derivativeLipschitzBound}\,R\,\|b-a\|\,\|h\|
\end{aligned}
\]
Uniform bounds on these two constants across the component expressions \(q_l\) supply
the constants in the Taylor \(\Znorm{1}\) and \(\Znorm{2}\) bounds
(Section~\ref{sec:taylor-services}).

\subsection{Finite/tail operators}\label{sec:finite-tail}

The polynomial-zero and Taylor families specify \(A\) and \(A^\dagger\) by a finite
block and a diagonal tail. With the truncation \(N\), each is an operator
\(B\) of the form
\[
 B=\begin{pmatrix}B^{(N)}&0\\0&\Lambda\end{pmatrix},\qquad
 (Bh)_n=\begin{cases}\displaystyle\sum_{k=0}^{N}B^{(N)}_{nk}h_k,&0\le n\le N,\\\Lambda_nh_n,&n>N.\end{cases}
\]
For systems, \(h=(h_l)_{0\le l<L}\), the finite block consists of \(L\times L\)
blocks \(B^{(N)}_{lj}\), each an \((N+1)\times(N+1)\) matrix on the modes
\(0,\ldots,N\) of the component \(h_j\), and \(\Lambda_n\) is diagonal on the
components. For \(B=A^\dagger\) and \(B=A\) the tail \(\Lambda_n\) is written
\(A^\dagger_{\mathrm{tail}}(n)\) and \(A_{\mathrm{tail}}(n)\). The tail of
\(A^\dagger\) is the diagonal of \(DF(\abar)\). For
\(n>N\),
\[
\begin{aligned}
 A^\dagger_{\mathrm{tail}}(n)&=2\abar_0,
 &A_{\mathrm{tail}}(n)&=(2\abar_0)^{-1}
 &&\text{(square root, }\abar_0\ne0),\\
 A^\dagger_{\mathrm{tail}}(n)&=n,
 &A_{\mathrm{tail}}(n)&=1/n
 &&\text{(Taylor IVP, in each component)},
\end{aligned}
\]
so that \(A_{\mathrm{tail}}(n)A^\dagger_{\mathrm{tail}}(n)=1\) for \(n>N\).
In both families the finite block of \(A^\dagger\) holds the entries of
\(DF(\abar)\) on the modes \(0,\ldots,N\) (for the Taylor family by the check
of Section~\ref{sec:taylor-services}), and since the derivative is lower
triangular, \(A^\dagger\) and \(DF(\abar)\) agree on the rows \(n\le N\). For a bounded operator \(B\), such as \(A\), the Lean data are
\begin{lstlisting}
structure SystemBlockDiagData (L N : ℕ) where
  finBlock : FiniteBlockMatrix L N
  tailDiag : Fin L → ℕ → ℝ
  tailBound : ℝ
  tailBound_spec : ∀ l n, N < n → |tailDiag l n| ≤ tailBound
\end{lstlisting}
where \lean{FiniteBlockMatrix}~\(L\)~\(N\) holds the blocks \(B^{(N)}_{lj}\),
\lean{tailDiag} the tail \(\Lambda\), and \lean{tailBound} a uniform bound on the absolute values of the tail
entries, written \(\tau_A\) for \(B=A\) (\(1/(N+1)\) for the Taylor initial value problems). The bounded
operator on
\((\ell^1_\nu)^L\), normed by
\(\|a\|=\max_l\|a_l\|_{1,\nu}\), is \lean{SystemBlockDiagData.toCLM}.

\paragraph{Assembly bound.} Equation~\eqref{eq:column-principle} bounds the operator norm
by the blocks and the tail bound. For every \(L\ge1\), \(N\) and \(\nu>0\),
\[
 \|A\|\le\max_{0\le l<L}\sum_{j=0}^{L-1}\bigl\|A^{(N)}_{lj}\bigr\|_{1,\nu}+\tau_A,
 \qquad
 \|H\|_{1,\nu}=\max_{0\le j\le N}\frac{\|H\delta_j\|_{1,\nu}}{\nu^j}
\]
where \(\|H\|_{1,\nu}\) is the finite weighted matrix norm of
Section~\ref{sec:free-banach}. For \(L=1\) a sharper bound is
\(\max\{\|A^{(N)}\|_{1,\nu},\tau_A\}\), the norm itself when \(\tau_A\) is
the supremum of the tail.

\paragraph{Composition with unbounded data.} For Taylor initial value
problems the tail of \(A^\dagger\) acts by \(\delta_n\mapsto n\delta_n\),
with \(\|n\delta_n\|_{1,\nu}/\nu^n=n\), so \(A^\dagger\) is block data and
defines no bounded operator. For a bounded \(A\) and data \(A^\dagger\) with
\(A_{\mathrm{tail}}A^\dagger_{\mathrm{tail}}=1\) on the tail,
\lean{SystemBlockDiagData.composedApprox} forms the data
\[
 (AA^\dagger)^{(N)}_{lj}=\sum_{m=0}^{L-1}A^{(N)}_{lm}(A^\dagger)^{(N)}_{mj},\qquad
 (AA^\dagger)_{\mathrm{tail}}=1,
\]
and \lean{toCLM} produces the bounded operator \(AA^\dagger\). No bound on
\(A^\dagger\) enters. With \(A_{\mathrm{tail}}(n)=1/n\) the Taylor
coefficient map~\eqref{eq:ivp-coeff} becomes
\[
 G(a)_{l,n+1}=\tfrac{1}{n+1}F(a)_{l,n+1}=a_{l,n+1}-\tfrac{\varphi_l(a)_n}{n+1}\qquad(n+1>N),
\]
and \(A(F(a))\in(\ell^1_\nu)^L\) for every \(a\), so that \(G=A\circ F\) is a map
\((\ell^1_\nu)^L\to(\ell^1_\nu)^L\) (\lean{ivpMap}).

\paragraph{Defect and injectivity.} The defect \(I-AA^\dagger\) has the
finite block \(I-(AA^\dagger)^{(N)}\) and a zero tail, so with
\(\|D\|_{\mathrm{sys}}=\max_l\sum_j\|D_{lj}\|_{1,\nu}\) for a finite block,
\[
 \bigl\|I-(AA^\dagger)^{(N)}\bigr\|_{\mathrm{sys}}\le\Znorm{0}
 \quad\Longrightarrow\quad\|I-AA^\dagger\|\le\Znorm{0}
\]
and the Taylor certificates check the
left side over \(\Q\) (Section~\ref{sec:taylor-ivp}). A strict bound below \(1\)
also gives injectivity. For matrices on the modes \(0,\ldots,N\),
\[
 \|I-M\|_{1,\nu}<1\ \Longrightarrow\ M\text{ is a unit},\qquad
 \bigl\|I-(AA^\dagger)^{(N)}\bigr\|_{\mathrm{sys}}<1\ \Longrightarrow\ A^{(N)}\text{ is a unit}
\]
With nonzero tail entries the block-diagonal \(A\) is
then injective, in the scalar case as a bounded operator and for systems on
raw coefficient families. The polynomial zero uses this as the injectivity premise of
Theorem~\ref{thm:radii-poly}, and the Taylor family uses it to pass from
\(G(\atilde)=0\) to \(F(\atilde)=0\) (Section~\ref{sec:taylor-ivp}).

\subsection{Scalar reductions}\label{sec:four-bounds}

For the polynomial zero, \(F(a)=P(a)-c\) on \(\ell^1_\nu\) with \(P\)
univariate and \(c\) finitely supported, each of the four bounds reduces to
a finite problem of its own kind.

\emph{Residual \(\Ynorm\).} Since \(\abar\) and \(c\) have finite support, so do
\(P(\abar)-c\) and \(A(P(\abar)-c)\). The weighted
norm of a sequence \(a\) with \(a_n=0\) for \(n>M\) is the finite sum
\(\sum_{n\le M}|a_n|\nu^n\).

\emph{Defect \(\Znorm{0}\).} Since \(AA^\dagger\) has identity tail, the columns of
\(I-AA^\dagger\) with index \(j>N\) vanish, and~\eqref{eq:column-principle} gives
\[
 \|I-AA^\dagger\|\le\max_{0\le j\le N}\frac{\|\bigl(I^{(N)}-A^{(N)}(A^\dagger)^{(N)}\bigr)\delta_j\|_{1,\nu}}{\nu^j}
\]

\emph{Derivative defect \(\Znorm{1}\).} Since \(DF(\abar)\) is lower triangular and
\((A^\dagger)^{(N)}=DF(\abar)^{(N)}\), the difference
\(\Delta=A^\dagger-DF(\abar)\) vanishes on the rows \(n\le N\), and on the
tail it is the difference of two left actions. In general, for bounded
\(\Delta,\Delta'\) with \((\Delta h)_n=0\) for \(n\le N\) and
\((\Delta h)_n=(\Delta'h)_n\) for \(n>N\),
\[
 \|A\circ\Delta\|\le\tau_A\,\|\Delta'\|
\]
For the
univariate polynomial \(P\) of the family and \(\Delta=A^\dagger-DF(\abar)\),
\(\Delta'\) is the left action of the element
\(-(P'(\abar)-P'(\abar)_0\cdot1)\), so one may take
\(\Znorm{1}=|A_{\mathrm{tail}}|\,\|P'(\abar)-P'(\abar)_0\cdot1\|_{1,\nu}\), a
finite sum over the support of \(P'(\abar)\).

\emph{Derivative variation \(\Znorm{2}\).} On a normed commutative ring
\(\mathcal A\) with a normed \(\R\)-algebra structure, write \(M_b\) for
the left action of the element \(b\). If
\(DF(x)=M_{g(x)}+\Psi\) for a fixed bounded linear map \(\Psi\), and if there is a
constant \(\kappa\ge0\) such that
\(\|g(x)-g(\xbar)\|\le\kappa\|x-\xbar\|\) on \(\closedball{r_0}{\xbar}\),
then
\[
 \|A\circ[DF(x)-DF(\xbar)]\|\le\kappa\,\|A\|\,r_0
 \qquad(x\in\closedball{r_0}{\xbar})
\]
For \(F(x)=P(x)-c\) on \(\ell^1_\nu\), \(DF(x)=M_{P'(x)}\), and with \(\|\abar\|_{1,\nu}+r_0\le R\) the family bounds
\(\|A\|\) by the assembly bound of Section~\ref{sec:finite-tail} and, for an
expression \(p\) that represents \(P\), takes
\(\kappa=p.\lean{derivativeLipschitzBound}\,R\). \lean{Example77} uses these four reductions
through its family theorem (Section~\ref{sec:sqrt}).

\subsection{Taylor estimates}\label{sec:taylor-services}

The certified map of the Taylor family is \(G=A\circ F\) of
Section~\ref{sec:fun-to-seq}, with the identity as preconditioner and
\(AA^\dagger\) as approximate derivative, so~\eqref{eq:Z1} reads
\(\|AA^\dagger-DG(\abar)\|\le\Znorm{1}\) and~\eqref{eq:Z2} reads
\(\|DG(c)-DG(\abar)\|\le\Znorm{2}r_0\) for \(c\in\closedball{r_0}{\abar}\).
Fix \(L\ge1\), a truncation \(N\) and a weight \(\nu>0\). The estimates
below apply to differentiable maps
\(\varphi\colon(\ell^1_\nu)^L\to(\ell^1_\nu)^L\) under the indicated bounds and
coefficient identities. For the expressions \(q_l\) of
Section~\ref{sec:fun-to-seq}, Section~\ref{sec:poly-derivatives} supplies the
derivative and its bounds.

\paragraph{Finite Jacobian.} The certificate stores the finite block of
\(DF(\abar)\) as a rational matrix, and the Jacobian identity compares it,
entry by entry over \(\Q\), with the derivative of~\eqref{eq:ivp-coeff}. Let
\(c\) denote the rational arrays representing \(\abar\). The rows with mode
zero are those of the identity. For \(n+1\le N\),
\[
 DF(\abar)_{(l,n+1),(j,k)}=(n+1)\,\delta_{lj}\,\delta_{n+1,k}
 -\begin{cases}
 (\partial_jq_l).\lean{evalCoeff}\,c\,(n-k),&k\le n,\\
 0,&k>n.
 \end{cases}
\]
The checked entries imply that \(AA^\dagger\) and \(DG(\abar)\) agree on the rows \(n\le N\), so
\(\bigl((AA^\dagger-DG(\abar))h\bigr)_{l,n}=0\) for every \(h\), \(l\) and
\(n\le N\).

\paragraph{The \(\Znorm{1}\) bound.} Let \(S\colon\ell^1_\nu\to\ell^1_\nu\)
be the shift \((Sx)_0=0\), \((Sx)_{n+1}=x_n/(n+1)\). On the
tail, \(AA^\dagger\) is the identity and \(G(a)_{l,n+1}=a_{l,n+1}-\varphi_l(a)_n/(n+1)\),
so the defect of the approximate derivative is \(S\) applied to the
derivative of \(\varphi\). If, for every \(h\) and \(l\),
\[
\begin{aligned}
 \bigl((AA^\dagger-DG(\abar))h\bigr)_{l,n}&=0&&(n\le N),\\
 \bigl((AA^\dagger-DG(\abar))h\bigr)_{l,n}&=\bigl(S\,D\varphi_l(\abar)h\bigr)_n&&(n>N),\\
 \textstyle\sum_{n>N}\bigl|\bigl(S\,D\varphi_l(\abar)h\bigr)_n\bigr|\,\nu^n&\le C\,\|h\|,
\end{aligned}
\]
with \(0\le C\le\Znorm{1}\), then \(\|AA^\dagger-DG(\abar)\|\le\Znorm{1}\)
(\lean{ivp\_Z\(_1\)\_le\_exact}). Once the finite rows cancel, the tail sums
are the component norms. Two ways supply \(C\). The
first bounds the tail of \(S\) by \(\nu/(N+1)\). From
\(\|D\varphi_l(\abar)h\|_{1,\nu}\le K\|h\|\) it takes \(C=\nu K/(N+1)\)
(\lean{ivp\_Z\(_1\)\_le}, cf.\ the Lorenz bound of
\cite[\S4.3]{kaliesKMJ2018}). The Lorenz certificate \lean{Example83}
(Section~\ref{sec:taylor-ivp}) uses it with \(K\) from
\lean{derivativeBound}. The second uses the columns of each scalar tail operator. When
\(D\varphi_l(\abar)h=\sum_mb_{lm}\star h_m\) with \(b_{lm}\in\ell^1_\nu\), the
tail sum is at most \(\sum_m\|W_{b_{lm}}\|\,\|h\|\) for
\(W_b=\pi_{>N}\circ S\circ M_b\), where \(\pi_{>N}\) sets the modes \(n\le N\) to zero. If \(b_k=0\) for
\(k>d\), the weighted column norms of \(W_b\) are
\[
 \gamma_j:=\frac{\|W_b\delta_j\|_{1,\nu}}{\nu^j}
 =\sum_{\substack{0\le k\le d\\ j+k+1>N}}\frac{|b_k|\,\nu^{k+1}}{j+k+1},
\]
nonincreasing in \(j\) for \(j\ge N\), so \(\|W_b\|=\max_{0\le j\le N}\gamma_j\)
by~\eqref{eq:column-principle}. Rational upper bounds
\(\gamma_j\le C_b\) for \(j\le N\) give \(\|W_b\|\le C_b\), and the scalar certificate \lean{Example81} checks
them.

\paragraph{The \(\Znorm{2}\) bound.} Let \(c\in\closedball{r_0}{\abar}\), with \(G\) and the
components of \(\varphi\) differentiable at \(c\) and \(\abar\). For a
direction \(h\), \((DG(c)-DG(\abar))h=Aw\), where \(w_l\) is
\(-(D\varphi_l(c)-D\varphi_l(\abar))h\) shifted up by one mode, so that
\(\|w_l\|_{1,\nu}=\nu\,\|(D\varphi_l(c)-D\varphi_l(\abar))h\|_{1,\nu}\). If an \emph{active set} \(J\subseteq\{0,\ldots,L-1\}\) and constants
\(\kappa,\Znorm{2}\ge0\) satisfy the following conditions for every \(h\), \(j\)
and \(l\),
\[
\begin{aligned}
 \bigl(D\varphi_j(c)-D\varphi_j(\abar)\bigr)h&=0&&(j\notin J),\\
 \bigl\|\bigl(D\varphi_l(c)-D\varphi_l(\abar)\bigr)h\bigr\|_{1,\nu}&\le \kappa\,\|c-\abar\|\,\|h\|,\\
 \kappa\,\nu\Bigl(\sum_{j\in J}\bigl\|A^{(N)}_{lj}\bigr\|_{1,\nu}+\mathbf 1_{l\in J}\,\tau_A\Bigr)&\le\Znorm{2},
\end{aligned}
\]
then \(\|DG(c)-DG(\abar)\|\le\Znorm{2}\,r_0\) (\lean{ivp\_Z\(_2\)\_le}). The
constant uses the assembly bound restricted to the columns in \(J\).
For the expressions \(q_l\) the rational clause
takes \(\kappa=\max_lq_l.\lean{derivativeLipschitzBound}\)
at the radii \(\|\abar_i\|_{1,\nu}+r_0\) and \(J\) as the components \(j\) with
\(q_j.\lean{derivativeLipschitzBound}\ne0\) at those radii. The Lorenz
certificate has \(J=\{1,2\}\), its first component being linear
(Section~\ref{sec:taylor-ivp}).

\subsection{Certificate checks}\label{sec:certificate-checks}

The reductions leave finite rational premises in the certificate data. Casting
from \(\Q\) to \(\R\) preserves the finite identities and inequalities. The
preceding estimates turn them into bounds on the real norms. For the polynomial-zero family the
soundness theorem \lean{PolyZeroData.certified} turns the decidable
proposition \(d.\lean{Checks}\) on the stored data \(d\) into \(r_0>0\),
\eqref{eq:Y}--\eqref{eq:Z2}, the differentiability of \(F\), \(p(r_0)<0\) and
the injectivity of \(A\), when \(\nu\) and the coefficients of \(\abar\) and of the
forcing \(c\) are the real casts of the stored rationals. For the
polynomial zero, \lean{decide~+kernel} proves \(d.\lean{Checks}\), which
includes the radii inequality \(\Znorm{2}r_0^2-(1-\Znorm{0}-\Znorm{1})r_0+\Ynorm<0\),
and \(\|A\|\) is bounded by the
exact rational value \(\max\{\|A^{(N)}\|_{1,\nu},|A_{\mathrm{tail}}|\}\).
Section~\ref{sec:verification} describes how the examples check their remaining
identities and inequalities. The Taylor family
states its clauses on the rational record
\begin{lstlisting}
structure StdIVPDataQ (L N : ℕ) [NeZero L] where
  A_col : Fin L → Fin L → ℕ → Array ℚ
  DF_col : Fin L → Fin L → ℕ → Array ℚ
  abar_Q : Fin L → Array ℚ
  ν_q : ℚ
  habar_size : ∀ l, (abar_Q l).size = N + 1
\end{lstlisting}
whose arrays hold the columns of \(A^{(N)}\) and \((A^\dagger)^{(N)}\) and
the coefficients of \(\abar\). The checks use the rational weight \(\nu_q\). The resulting
estimates use its real cast \(\nu\).

\section{Examples}\label{sec:examples}

Each example certifies one zero with the estimates of
Section~\ref{sec:reductions} and then returns to the equation it came from,
along the path
\[
\begin{aligned}
 &\text{data}\ \longrightarrow\ \bigl(\text{four norm bounds and }p(r_0)<0\bigr)\\
 &\qquad\longrightarrow\ \text{zero of the certified map}\ \longrightarrow\ \text{realization}.
\end{aligned}
\]
For the polynomial zero, evaluation turns the coefficient zero into the
square-root branch. For the
Taylor initial value problems, the constant coefficients give the initial values
and termwise differentiation turns the remaining coefficient equations into the differential equation.

\subsection{Polynomial zero}\label{sec:sqrt}

The unknown is a power series \(x(t)=\sum_{n\ge0}a_nt^n\) with
\(x(t)^2=\lambda_0+t\), and the certified map acts on its coefficient
sequence \(a\in\ell^1_\nu\) (Section~\ref{sec:fun-to-seq}).

\paragraph{Data.}
In the polynomial-zero family, a univariate expression \(p\) represents the
polynomial \(P\), and the coefficient equation is \(P(a)=c\) for finitely
supported \(c\). Here \(p=X_0\cdot X_0\) represents \(P=X_0^2\), so
\(F(a)=a\star a-c\), with \(c=(\lambda_0,1,0,\ldots)\).
The record \lean{RadiiPolynomial.PolyZeroData} stores \(p\), the truncation
\(N\), and rational data for the weight \(\nu_q\), the coefficients of
\(\abar\) and \(c\), the finite columns of \(A\), the radius and the proposed
bounds. It computes \(A_{\mathrm{tail}}=(A^\dagger_{\mathrm{tail}})^{-1}\).
The record \lean{Example77.Cert.data} has \(\lambda_0=1/3\),
\(\nu=1/4\) and \(N=2\), with
\[
 \abar=\bigl(\tfrac{5774}{10000},\ \tfrac{8660}{10000},\ -\tfrac{6495}{10000}\bigr),\qquad
 A^{(2)}=\frac{1}{10000}\begin{pmatrix}
 8660&0&0\\
 -12990&8660&0\\
 29240&-12990&8660
 \end{pmatrix},\qquad
 r_0=\tfrac{996}{10000}.
\]

\paragraph{Bounds.}
The four reductions of Section~\ref{sec:four-bounds} give, for
\(a\in\closedball{r_0}{\abar}\),
\[
\begin{aligned}
 \|A(F(\abar))\|_{1,\nu}&\le\Ynorm=\tfrac{9}{500},\\
 \|I-AA^\dagger\|&\le\Znorm{0}=\tfrac{2}{1000},\\
 \|A\circ[A^\dagger-DF(\abar)]\|&\le\Znorm{1}=\tfrac{46}{100},\\
 \|A\circ[DF(a)-DF(\abar)]\|&\le\Znorm{2}\,r_0,\qquad\Znorm{2}=\tfrac{28}{10}.
\end{aligned}
\]
The value \(\Znorm{1}\) bounds
\(|A_{\mathrm{tail}}|\,\|P'(\abar)-P'(\abar)_0\cdot1\|_{1,\nu}\).
At \(R=\|\abar\|_{1,\nu}+r_0\), the expression gives
\(p.\lean{derivativeLipschitzBound}\,R=2\). With the sharper bound on
\(\|A\|\) for \(L=1\) from Section~\ref{sec:finite-tail}, the check for \(\Znorm{2}\) is
\[
 2\max\{\|A^{(2)}\|_{1,\nu},|A_{\mathrm{tail}}|\}
 =\tfrac{2747}{1000}\le\tfrac{28}{10}.
\]

\paragraph{Radii inequality and zero.}
\begin{theorem}[Sequence-space square root, \lean{Example77.Cert.main\_theorem}]\label{thm:example77}
There exists a unique \(\atilde \in \ell^1_\nu\) with
\(\|\atilde - \abar\|_{1,\nu} \leq 996/10000\) satisfying
\(\atilde \star \atilde = c\).
\end{theorem}

The proposition \lean{PolyZeroData.Checks} compares the four rational
computations with the proposed bounds and includes the inequality
\(\Znorm{2}r_0^2-(1-\Znorm{0}-\Znorm{1})\,r_0+\Ynorm<0\).
The proof \lean{Example77.Cert.}\allowbreak\lean{data\_checks} uses one
\lean{decide~+kernel}. The family theorem
\lean{RadiiPolynomial.}\allowbreak\lean{PolyZeroData.}\allowbreak\lean{existsUnique} applies
Theorem~\ref{thm:radii-poly} to \(F\) itself, with the injectivity of \(A\)
from Section~\ref{sec:finite-tail}, and returns the unique
\(\atilde\in\closedball{r_0}{\abar}\) with \(P(\atilde)=c\).

\paragraph{Realization.}
Put \(x(t)=\lean{eval}(\atilde,t)=\sum_n\atilde_nt^n\), which converges
absolutely for \(|t|\le\nu\) by Definition~\ref{def:l1Weighted}. Then
\[
 x(t)^2=\lean{eval}(\atilde\star\atilde,t)
       =\lean{eval}(c,t)=\lambda_0+t,
\]
where the first equality is the substitution identity of
Section~\ref{sec:free} for \(P=X_0^2\) at the real point \(t\).
For this polynomial it follows from Theorem~\ref{thm:mertens}.
The next two equalities use \(\atilde\star\atilde=c\) and
\(c=(\lambda_0,1,0,\ldots)\), respectively.

\begin{theorem}[Square-root identity, \lean{Example77.Cert.x\_analytic\_sq\_eq}]
\label{thm:example77-analytic}
Let \(x(t)=\sum_{n} \atilde_n t^n\). For \(t\in\R\) with
\(|t|\le\nu\), one has \(x(t)^2 = \lambda_0 + t\).
\end{theorem}

The series defines an analytic branch for \(|t|<\nu\). The coefficient ball
selects the positive square root at \(t=0\), since
\(x(0)=\atilde_0\ge\abar_0-r_0>0\).

\subsection{Taylor initial value problems}\label{sec:taylor-ivp}

In the Taylor family, the factor \(n+1\) in the coefficient
map~\eqref{eq:ivp-coeff} leads to the preconditioned map \(G\) of
Section~\ref{sec:finite-tail}. The scalar certificate \lean{Example81} and the
Lorenz certificate \lean{Example83} use the same family theorems, with
different estimates for \(\Znorm{1}\) from Section~\ref{sec:taylor-services}.

\paragraph{Data.}
The record \lean{StdIVPData}\,\(\nu\)\,\(L\)\,\(N\) contains the fields of
\lean{StdIVPDataQ} from Section~\ref{sec:certificate-checks} and the equality
\(\nu=(\nu_q:\R)\). The rational arrays give the coefficients of \(\abar\)
and the entries of \(A^{(N)}\) and \((A^\dagger)^{(N)}\).
For \(n>N\) the tails are \(A^\dagger_{\mathrm{tail}}(n)=n\) and
\(A_{\mathrm{tail}}(n)=1/n\) in every component.
The expressions \(q_l=\lean{f\_cpoly}_l\) give both
\(\varphi_l(a)=q_l.\lean{evalBanach}\,a\) for \(a\in(\ell^1_\nu)^L\)
and \(f_l(y)=q_l.\lean{evalBanach}\,y\) for \(y\in\R^L\).

The scalar certificate \lean{Example81} treats \(\dot x=x(x-1)\),
\(x(0)=1/2\), with \(L=1\), \(N=10\) and \(\nu=1\). Its \(A^{(N)}\) is the
exact rational inverse of the finite Jacobian \((A^\dagger)^{(N)}\), so
\(I-(AA^\dagger)^{(N)}=0\). For \(y=(y_0,y_1,y_2)\in\R^3\), the Lorenz field is
\[
 \begin{aligned}
 f_0(y)&=10(y_1-y_0),\\
 f_1(y)&=28y_0-y_1-y_0y_2,\\
 f_2(y)&=-(8/3)y_2+y_0y_1,
 \end{aligned}
\]
with the initial condition \(x(0)=(1,0,0)\), truncation \(N=30\) and weight
\(\nu=3/20\). An external computation provides the rational coefficients
in \((\Q^{31})^3\) and two matrices of blocks in
\((\Q^{31\times31})^{3\times3}\). The certificate checks the stored Jacobian
against the partial derivatives computed from \(q_l\), as in
Section~\ref{sec:taylor-services}. For the Lorenz estimate of \(\Znorm{1}\),
the component norms of \(\abar\) satisfy
\[
 \mathbf R=(20,26,11),\qquad
 \|\abar_l\|_{1,\nu}\le\mathbf R_l .
\]

\paragraph{Bounds.}
For \(c\in\closedball{r_0}{\abar}\), the four norm bounds are
\[
\begin{aligned}
 \|G(\abar)\|&\le\Ynorm,\\
 \bigl\|I-(AA^\dagger)^{(N)}\bigr\|_{\mathrm{sys}}&\le\Znorm{0},\\
 \|AA^\dagger-DG(\abar)\|&\le\Znorm{1},\\
 \|DG(c)-DG(\abar)\|&\le\Znorm{2}\,r_0,
\end{aligned}
\]
with the following certified values
\[
\begin{array}{@{}lccccc@{}}
 & \Ynorm & \Znorm{0} & \Znorm{1} & \Znorm{2} & r_0\\[2pt]
 \lean{Example81} & \tfrac{3}{10^{6}} & 0 & \tfrac{482763031}{9859449600} & \tfrac{110916221}{39916800} & \tfrac{1}{10}\\[4pt]
 \lean{Example83} & \tfrac{44}{10^{5}} & \tfrac{37}{10^{15}} & \tfrac{30}{100} & \tfrac{49}{10} & \tfrac{64}{10^{5}}
\end{array}
\]
Both certificates prove the finite inequality for \(\Znorm{0}\) by a rational
check with \lean{decide~+kernel}. For \(\Znorm{1}\), the scalar certificate uses the
largest weighted column norm of the tail operator in Section~\ref{sec:taylor-services}. The
Lorenz certificate takes \(C=\nu K/(N+1)\) with \(K=62\), since
\(q_l.\lean{derivativeBound}\,\mathbf R\le62=\Znorm{1}(N+1)/\nu\)
for every component. For \(\Znorm{2}\), the values of \lean{derivativeLipschitzBound} for the three
Lorenz expressions are \(0,2,2\), respectively, giving \(J=\{1,2\}\). The scalar
expression \(X_0\cdot X_0-X_0\) gives \(J=\{0\}\).

\paragraph{Radii inequality and zero.}
The inequality \(\Znorm{2}r_0^2-(1-\Znorm{0}-\Znorm{1})\,r_0+\Ynorm<0\)
holds at the value of \(r_0\) in each row of the table. The family theorem
\lean{IVP.StdIVPData.}\allowbreak\lean{existsUnique\_}\allowbreak\lean{of\_compPoly}
applies Theorem~\ref{thm:radii-poly} to \(G\), with \((F,A,A^\dagger)\)
replaced by \((G,I,AA^\dagger)\), and gives the unique
\(\atilde\in\closedball{r_0}{\abar}\) with \(G(\atilde)=0\).
Since \(r_0>0\) and the norm bounds are nonnegative, \(p(r_0)<0\) implies
\(\Znorm{0}<1\). With \(A_{\mathrm{tail}}(n)=1/n\ne0\) for \(n>N\),
Section~\ref{sec:finite-tail} gives injectivity of \(A\) on raw coefficient
families. Thus \(G(a)=0\) if and only if every coefficient of \(F(a)\)
vanishes. The family theorem
\lean{IVP.StdIVPData.}\allowbreak\lean{existsUnique\_}\allowbreak\lean{ivpCoeffs\_of\_compPoly}
gives a unique zero of the coefficient equations in the same ball.
Its scalar instance is
\lean{Example81.}\allowbreak\lean{Cert.}\allowbreak\lean{ivp\_main\_theorem},
and its Lorenz instance is the following theorem.

\begin{theorem}[Lorenz coefficient solution,
  \lean{Example83.Cert.ivp\_main\_theorem}]
\label{thm:lorenz-coeff}
There exists a unique
\(\atilde\in\overline B_{r_0}(\abar)\) such that
\[
  F(\atilde)_{l,n}=0
  \qquad\text{for }l=0,1,2\text{ and every }n\in\N.
\]
\end{theorem}

\paragraph{Realization.}
Choose \(\atilde\) to be the certified zero of the coefficient
map~\eqref{eq:ivp-coeff} and evaluate componentwise,
\(x(t)_l=\lean{eval}(\atilde_l,t)\). The row \(F(\atilde)_{l,0}=0\) gives
\(x(0)_l=(x_0)_l\), and
\begin{equation}\label{eq:ode-chain}
\begin{aligned}
 \dot x_l(t)
 &=\sum_{n\ge0}(n+1)\atilde_{l,n+1}t^n
  =\sum_{n\ge0}\varphi_l(\atilde)_n t^n\\
 &=\lean{eval}(\varphi_l(\atilde),t)
  =f_l(x(t)),\qquad |t|<\nu,
\end{aligned}
\end{equation}
by termwise differentiation (\lean{hasDerivAt\_eval}), the equations
\(F(\atilde)_{l,n+1}=0\), and polynomial substitution at \(t\)
(\lean{eval\_evalBanach}, Section~\ref{sec:free}). Together with
\(x(0)=x_0\), this gives \lean{IVP.}\allowbreak\lean{solves\_ODE\_of\_F\_zero}.
For \(|t|\le\nu\), evaluation bounds \(|x(t)_l|\) by
\(\|\atilde_l\|_{1,\nu}\), and \(\atilde\in\closedball{r_0}{\abar}\) gives
\[
 \|x(t)\|_\infty\le\|\atilde\| \le\|\abar\|+r_0=R,
\]
where \(R\) is the library's \lean{R\_traj}. On the ball
\(\closedball{R}{0}\subseteq\R^L\), the expressions \(q_l\) give
\[
\begin{aligned}
 K_f&:=\sum_{l=0}^{L-1}\bigl|q_l.\lean{lipschitzBound}\,(R,\dots,R)\bigr|,\\
 \|f(y)-f(y')\|_\infty&\le K_f\|y-y'\|_\infty
 \qquad(y,y'\in\closedball{R}{0}),
\end{aligned}
\]
so every differentiable solution on \((-\nu,\nu)\) with initial value
\(x_0\) and trajectory in \(\closedball{R}{0}\) agrees there with \(x\)
(\lean{IVP.}\allowbreak\lean{analytic\_}\allowbreak\lean{solution\_unique}).
The family theorem
\lean{IVP.StdIVPData.}\allowbreak\lean{analytic\_existsUnique\_}\allowbreak\lean{of\_compPoly\_}\allowbreak\lean{of\_radii}
gives existence and uniqueness on the interval from the four bounds at
\(r_0>0\) and the inequality \(p(r_0)<0\). Its scalar instance is
\lean{Example81.Cert.}\allowbreak\lean{analytic\_existsUnique} on \((-1,1)\).
Its Lorenz instance is the following theorem.

\begin{theorem}[Lorenz IVP, \lean{Example83.Cert.analytic\_existsUnique}]
\label{thm:lorenz-analytic}
The Lorenz initial value problem with \(x(0)=(1,0,0)\) has a solution
on \((-\nu,\nu)\), with trajectory in \(\closedball{R}{0}\), given componentwise by the power series
\(x(t)_l=\sum_{n=0}^{\infty}(\atilde_l)_n t^n\). Every solution of this initial value problem with trajectory in \(\closedball{R}{0}\) agrees with it on \((-\nu,\nu)\).
\end{theorem}

The proof takes \(u=x\), the series solution constructed above.
In the namespace \lean{Example83.Cert} the Lean statement reads
\begin{lstlisting}
theorem analytic_existsUnique :
    ∃ u : ℝ → Fin L → ℝ,
      IVP.IsBallSolution (ν := ν_val) f_cpoly x₀ R_traj u ∧
      ∀ v : ℝ → Fin L → ℝ,
        IVP.IsBallSolution (ν := ν_val) f_cpoly x₀ R_traj v →
        Set.EqOn v u (Set.Ioo (-(ν_val : ℝ)) ν_val)
\end{lstlisting}
with \(\lean{L}=3\), \(\lean{\(\nu\)\_val}=\nu\) and \(\lean{R\_traj}=R\). For a function \(y\colon\R\to\R^L\), the
predicate \lean{IVP.IsBallSolution} requires \(y\) to be differentiable on
\((-\nu,\nu)\), with
\[
 y(0)=x_0,\qquad \|y(t)\|_\infty\le R,\qquad \dot y(t)=f(y(t))\qquad(t\in(-\nu,\nu)),
\]
and no analyticity condition on \(y\). The series solution \(x\) is also
analytic on \((-\nu,\nu)\), by
\lean{Example83.Cert.}\allowbreak\lean{lorenz\_analytic\_existsUnique}.

\section{Verification}\label{sec:verification}

The certificates in Section~\ref{sec:examples} prove real norm bounds
from finite calculations. The assumptions of the resulting theorems depend
on how Lean checks those calculations.

With \lean{decide~+kernel}, Lean's kernel evaluates the decision procedure
for the proposition. With \lean{native\_decide}, compiled code evaluates
the decision procedure, and Lean records the successful Boolean check as
an auxiliary axiom. The kernel checks the proof using that axiom, so the
theorem also relies on Lean's compiler and runtime.\footnote{See the Lean reference manual, \emph{Validating a Lean
Proof}, \url{https://lean-lang.org/doc/reference/latest/ValidatingProofs/}.}
For the finite weighted sums in the Taylor certificates,
LeanCert~\cite{leancert} combines interval bounds for the individual terms
with a native check on their sum to prove the real inequality.

Both Taylor certificates check the finite Jacobian identity over \(\Q\).
The scalar certificate uses \lean{decide~+kernel}, and the Lorenz
certificate uses \lean{native\_decide}. The two radii inequalities are
stated over \(\R\), and neither proof requires a native axiom.

For the \(\Znorm{0}\) bound, both Taylor certificates use
\lean{decide~+kernel} with grid spacing \(10^{-60}\). The scalar
certificate computes the defect columns exactly before rounding the
absolute values of their entries up to the grid, so exact zeros remain
zero. The Lorenz certificate uses integer arithmetic on the grid to
obtain an upper bound for the finite defect norm.

\paragraph{Axioms.}\label{sec:axiom-profiles}
At revision \lean{\LibCommit}, \lean{\#print~axioms} lists the axioms
on which a theorem and its dependencies rely. Theorems~\ref{thm:radii-poly}
and~\ref{thm:radii-local} depend only on Lean's standard axioms
\lean{propext}, \lean{Classical.choice} and \lean{Quot.sound}. The same
holds for the estimates of Sections~\ref{sec:four-bounds}
and~\ref{sec:taylor-services}, the soundness and family theorems used in
Section~\ref{sec:examples}, and the polynomial-zero certificate with
Theorems~\ref{thm:example77} and~\ref{thm:example77-analytic}.

For the scalar Taylor example, the theorem giving the coefficient zero
and the theorem giving the solution on \((-1,1)\) depend on the same
additional axiom, from the check on the residual bound. The corresponding
Lorenz results, Theorems~\ref{thm:lorenz-coeff} and~\ref{thm:lorenz-analytic},
and the theorem asserting analyticity of the series solution depend on the
same seven additional axioms. Three come from the component bounds on the
residual, one from the finite Jacobian identity, and three from the bounds
\(\|\abar_l\|_{1,\nu}\le\mathbf R_l\).
None of these theorems depends on \lean{sorryAx}. The library source
contains no explicit axiom declarations.
\section{Discussion}\label{sec:discussion}

\paragraph{What the framework delivers.}
The operator estimates of Section~\ref{sec:reductions} give the norm bounds
used by Theorem~\ref{thm:radii-poly} to certify a zero in coefficient space.
Evaluation then gives a function. For the square root, evaluation preserves
the Cauchy product, so the coefficient equation gives \(x(t)^2=\lambda_0+t\).
For the Taylor initial value problems, injectivity of \(A\) first gives
\(F(\atilde)=0\) from \(G(\atilde)=0\). Termwise differentiation and the
substitution identity of Section~\ref{sec:free} then give
\(\dot x(t)=f(x(t))\), with the initial value given by the constant coefficients.

\paragraph{Reuse.}
For a new polynomial in either family, the evaluation and derivative
identities of Section~\ref{sec:free} need no new proof. The certificate
supplies the numerical data and proves the bounds by the finite calculations
of Section~\ref{sec:verification}. In Section~\ref{sec:taylor-ivp}, for
example, \lean{Example81} and \lean{Example83} use different estimates for
\(\Znorm{1}\) and apply the same Taylor family theorem. A new family can
use Theorem~\ref{thm:radii-poly} once its coefficient map and the required
bounds and hypotheses are established. To obtain a solution as a function,
it also needs a realization theorem.

\paragraph{Further applications.}
At revision \lean{\LibCommit}, \lean{StdFinDimData} certifies zeros of
polynomial systems with rational coefficients on \(\R^L\) when the
certificate checks pass. With \(L=1\), the library answers \textbf{P1} for
\(x^2-2\) by certifying a unique zero in the ball of radius \(3/20\)
about \(13/10\).
The library also certifies an analytic solution of an initial value problem
on \([-1,1]\) using Chebyshev coefficients, with convolution on \(\Z\).
Using weighted \(\ell^2\) Fourier coefficients, it certifies real analytic
periodic solutions of the forced quadratic Helmholtz--Duffing and Riccati
equations.

For boundary value problems and connecting orbits, the coefficient map must
also encode the boundary conditions~\cite{lessardReinhardt2014,vdbMJLM2011}.
Invariant manifolds can be computed from equations for their
parameterizations~\cite{kaliesKMJ2018}. Extending the library to certify
such boundary value solutions, connecting orbits and invariant manifolds
is left to future work.

\paragraph{Specification by universal properties.}
A bounded linear map out of \(\ell^1_\nu\) is determined by its values
on the atoms. When the assigned values satisfy the growth condition of
Section~\ref{sec:free-banach} and preserve the unit and products, their
linear lift is an algebra homomorphism. Evaluation at \(t\) is one such
map. The square-root proof can therefore invoke multiplicativity without
proving the Cauchy product formula again. These facts are classical.
The advantage for the library is that the same map and its laws can be
used for different equations.

Under the weight bounds of Appendix~\ref{app:lp}, the same algebra lift
sends the atoms of weighted \(\ell^1\) to shifts on weighted \(\ell^p\).
It gives an action by convolution without requiring a product of two
\(\ell^p\) elements.

\paragraph{Provenance.}
Mathlib supplies the polynomial algebra and the analysis used here. This
project contributed to Mathlib's discrete convolution interface and
supplied the lemma \lean{fderiv\_apply}. It also contributed the tactic
\lean{finsum\_bound} to LeanCert~\cite{leancert}, which supplies the
interval arithmetic.

\subsection{Disclosure of the use of artificial intelligence}

By now AI, or more concretely agentic systems, is an everyday instrument of formalization. Itemizing its uses would say little about its role in the research process. The models used most were Anthropic's Fable 5.1 and OpenAI's GPT-6 Astra, each run through its vendor's official harness. Handing everything to an agent is fine if all one wants is a type-checked statement. It goes badly wrong if what one wants is a well-structured API. At the time of writing, API design remains the hard part of autoformalization.  

\appendix
\section{Categorical formulation}\label{app:categorical}

For an additive monoid \(M\) and a positive submultiplicative weight
\(\omega\), the two lifts of Section~\ref{sec:algebras} describe maps out of
the same Banach space \(\ell^1_\omega(M)\). The linear lift makes this space
the coproduct of real lines in which a scalar \(c\) at index \(m\) has norm
\(|c|\omega_m\). A linear contraction into \(Y\) is determined by the images
\(y_m\) of the atoms, with \(\|y_m\|\le\omega_m\). When \(Y\) is a unital
real Banach algebra, the algebra lift gives a contractive unital algebra
homomorphism exactly when these images also satisfy \(y_0=1_Y\) and
\(y_{i+j}=y_i y_j\). We use the terminology of Mac Lane~\cite{maclane1998}
for universal arrows and adjunctions.

\paragraph{One object.}
Let \(U\colon\mathbf{Ban}_1\to\mathbf{wSet}_1\) send a Banach space to its
points weighted by the norm. Here \(\mathbf{wSet}_1\) is the category of
weighted sets \((X,w)\), \(w\colon X\to[0,\infty)\), whose arrows are the
weight decreasing maps, \(w'(f(x))\le w(x)\), and \(\mathbf{Ban}_1\) that of
real Banach spaces and linear contractions. An arrow
\(y\colon(M,\omega)\to U(Y)\) is a family of vectors with
\(\|y_m\|\le\omega_m\), and the atoms form the arrow
\(\eta\colon(M,\omega)\to U(\ell^1_\omega(M))\). The general form of the
lift displayed in Section~\ref{sec:free-banach} gives the natural bijection
\[
 \operatorname{Hom}_{\mathbf{Ban}_1}\bigl(\ell^1_\omega(M),-\bigr)
 \overset{\sim}{\Longrightarrow}
 \operatorname{Hom}_{\mathbf{wSet}_1}\bigl((M,\omega),U(-)\bigr),
 \qquad \psi\longmapsto(\psi(\delta_m))_{m\in M},
\]
so \(\ell^1_\omega(M)\) is the value at \((M,\omega)\) of a left adjoint
\(F\dashv U\). A weighted set is the coproduct of its points, and \(F\)
preserves coproducts. Thus \(F(M,\omega)\) is the coproduct in
\(\mathbf{Ban}_1\) of the lines \(\R_m\) with norm \(|c|\omega_m\) and
injections \(c\mapsto c\,\delta_m\). The map \(c\mapsto c\,y_m\) on
\(\R_m\) is a contraction precisely when \(\|y_m\|\le\omega_m\).
The displayed bijection therefore gives the universal property of this
coproduct.

Among real normed spaces and linear contractions, the coproduct is the
direct sum \(\bigoplus_m\R_m\) with norm \(\sum_m|a_m|\omega_m\).
Completion is left adjoint to the inclusion of real Banach spaces into
real normed spaces with those same arrows, so completing this direct sum
gives the Banach space coproduct. The direct sum is the underlying vector
space of the monoid algebra \(\R[M]\). Its multiplication is the bilinear
extension of \(\delta_i\star\delta_j=\delta_{i+j}\).
Submultiplicativity gives \(\|a\star b\|\le\|a\|\,\|b\|\) on the direct
sum, so multiplication extends to the completion. Conversely, evaluating
this inequality on two atoms gives \(\omega_{i+j}\le\omega_i\omega_j\).

In Lean the two lifts use the identifications \(\theta_m\) of
Section~\ref{sec:weighted-space}. Each \(\theta_m(1)\) is a generator
of its fiber and maps to \(\delta_m\) under the coproduct injection.
Section~\ref{sec:monoid-algebra} uses these same atoms for the
multiplication of the completed monoid algebra.

\paragraph{The bridge.}
Write \(\operatorname{Mult}(M)\) for \(M\) in multiplicative notation.
In \(\mathbf{wMon}_1\) an object is a monoid with a weight in
\([0,\infty)\), and an arrow is a weight decreasing monoid homomorphism.
In \(\mathbf{BanAlg}_1\) the objects are unital real Banach algebras, with
no condition on \(\|1\|\), and the arrows are contractive unital real
algebra homomorphisms. The functor
\(U_1\colon\mathbf{BanAlg}_1\to\mathbf{wMon}_1\) keeps the multiplicative
monoid and takes the norm as its weight. For our \((M,\omega)\), the map
\(\eta(m)=\delta_m\) is an arrow of \(\mathbf{wMon}_1\) because
\(\delta_0=1\) and \(\delta_i\star\delta_j=\delta_{i+j}\).
The algebra lift gives the universal arrow
\[
\begin{tikzcd}[column sep=large,row sep=large,ampersand replacement=\&]
 (\operatorname{Mult}(M),\omega)\arrow[r,"\eta"]\arrow[dr,"y"']
 \&U_1(\ell^1_\omega(M))\arrow[d,"U_1(\hat y)",dashed]
 \&\ell^1_\omega(M)\arrow[d,"{\exists!\,\hat y}",dashed]\\
 \&U_1(Y)\&Y
\end{tikzcd}
\]
with the natural bijection
\[
 \operatorname{Hom}_{\mathbf{BanAlg}_1}\bigl(\ell^1_\omega(M),-\bigr)
 \overset{\sim}{\Longrightarrow}
 \operatorname{Hom}_{\mathbf{wMon}_1}\bigl((\operatorname{Mult}(M),\omega),U_1(-)\bigr),
 \qquad \psi\longmapsto(\psi(\delta_m))_{m\in M}.
\]
This bijection restricts the linear one to families satisfying
\(y_0=1_Y\) and \(y_{i+j}=y_i y_j\). For the linear lift \(\hat y\),
the bilinear maps \((a,b)\mapsto\hat y(a\star b)\) and
\((a,b)\mapsto\hat y(a)\hat y(b)\) agree on pairs of atoms. By bilinearity
they agree on pairs of finitely supported families, and continuity extends
the equality to the completion. The library proves multiplicativity by
applying uniqueness of the linear lift in each variable
(Section~\ref{sec:monoid-algebra}). Theorem~\ref{thm:mertens} uses the
resulting algebra homomorphism for evaluation.

For each fixed \(C\ge1\), a lift bounded by \(C\) is a contraction when
the source weight is \(C\omega\). This weight remains submultiplicative.
The algebra and its product are unchanged, and its norm is multiplied by
\(C\).

\paragraph{The instances.}
For \(M=\N\), the index monoid is free on \(1\).
A bounded unital real algebra homomorphism out of \(\ell^1_\nu\) is
determined by the image \(x\) of \(\delta_1\), since the polynomials
in \(\delta_1\) are dense. Section~\ref{sec:evaluations} constructs this
homomorphism from a bound on the powers of \(x\).
Let \(\mathbf{BanAlg}_\infty\) have the same objects as \(\mathbf{BanAlg}_1\)
and bounded unital real algebra homomorphisms as arrows.
With \(\Pow^1_\nu(Y)=\{x\in Y : \|x^n\|\le\nu^n\ \text{for all } n\}\) and
\(\Pow_\nu(Y)=\{x\in Y : \exists C\ge0\ \forall n,\ \|x^n\|\le C\,\nu^n\}\),
\[
 \operatorname{Hom}_{\mathbf{BanAlg}_1}(\ell^1_\nu,-)
 \overset{\sim}{\Longrightarrow}\Pow^1_\nu(-),\qquad
 \operatorname{Hom}_{\mathbf{BanAlg}_\infty}(\ell^1_\nu,-)
 \overset{\sim}{\Longrightarrow}\Pow_\nu(-),\qquad
 \psi\mapsto\psi(\delta_1),
\]
natural in \(Y\), with inverse \(x\mapsto\operatorname{ev}_x\).
For a commutative ring \(R\) and an index set \(\sigma\), the polynomial
algebra \(R[X_i:i\in\sigma]\) is the free commutative \(R\)-algebra on
\(\sigma\). Write \(\mathbf{CAlg}_R\) for the category of commutative
\(R\)-algebras and \(R\)-algebra homomorphisms, and
\(U_{\mathrm{Set}}\) for its forgetful functor to sets. Then
\[
 \operatorname{Hom}_{\mathbf{CAlg}_R}\bigl(R[X_i : i\in\sigma],-\bigr)
 \overset{\sim}{\Longrightarrow}
 \operatorname{Hom}_{\mathbf{Set}}\bigl(\sigma,U_{\mathrm{Set}}(-)\bigr),
 \qquad\psi\longmapsto(i\mapsto\psi(X_i)),
\]
with no norm condition on the images of the variables. For a commutative
Banach algebra \(Y\), the next diagram composes the two induced algebra
homomorphisms to give the substitution identity of Section~\ref{sec:intro}.

\paragraph{The composite.}
Take \(R=\R\) and \(\sigma=\{0,\dots,L-1\}\), and write
\(\R[X]=\R[X_0,\dots,X_{L-1}]\). Definition~\ref{def:aeval} at
\(S=\ell^1_\nu\) turns \(L\) sequences
\(a=(a_0,\dots,a_{L-1})\) into \(\operatorname{ev}_a\colon\R[X]\to\ell^1_\nu\).
The evaluation of Section~\ref{sec:evaluations} turns one element
\(x\in\Pow_\nu(Y)\) of a commutative complete normed \(\R\)-algebra \(Y\)
into \(\operatorname{ev}_x\colon\ell^1_\nu\to Y\). The \(L\) variables
\(X_i\) are not the one generator \(\delta_1\) of \(\ell^1_\nu\). The
datum \(a\) assigns to each variable \(X_i\) an arbitrary
element \(a_i\) of \(\ell^1_\nu\), and the datum \(x\) is the image of the
generator \(\delta_1\) in \(Y\). Since \(\ell^1_\nu\) and \(Y\)
are commutative, \(\operatorname{ev}_x\) is an arrow of \(\mathbf{CAlg}_\R\),
and Definition~\ref{def:aeval} at \(S=Y\) gives a third homomorphism
\(\operatorname{ev}_{\operatorname{ev}_x\circ a}\). The diagram below lives in
\(\mathbf{Set}\). Its two arrows out of \(\{0,\dots,L-1\}\) assign the
variables and the sequences. The other three are underlying maps of
\(\R\)-algebra homomorphisms. Definition~\ref{def:aeval} induces the
horizontal arrow from \(a\) and the lower diagonal arrow from
\(\operatorname{ev}_x\circ a\). The right vertical arrow is evaluation at
\(x\).
\[
\begin{tikzcd}[column sep=large,row sep=large,ampersand replacement=\&]
 \{0,\dots,L-1\}\arrow[d,"{i\mapsto X_i}"']\arrow[dr,"a"]
 \&\\
 U_{\mathrm{Set}}(\R[X])\arrow[r,"U_{\mathrm{Set}}(\operatorname{ev}_a)"]\arrow[dr,"U_{\mathrm{Set}}(\operatorname{ev}_{\operatorname{ev}_x\circ a})"']
 \&U_{\mathrm{Set}}(\ell^1_\nu)\arrow[d,"U_{\mathrm{Set}}(\operatorname{ev}_x)"]\\
 \&U_{\mathrm{Set}}(Y)
\end{tikzcd}
\]
For a polynomial vector field \(f=(f_0,\dots,f_{L-1})\) with
\(f_l\in\R[X]\), the lower triangle reads
\[
 \operatorname{ev}_x\bigl(\operatorname{ev}_a(f_l)\bigr)
 =\operatorname{ev}_{\operatorname{ev}_x\circ a}(f_l)
 =f_l\bigl(\operatorname{ev}_x(a_0),\dots,\operatorname{ev}_x(a_{L-1})\bigr)
 \qquad(0\le l<L).
\]
On the left, \(f_l\) is evaluated on the sequences and the
resulting series is summed at \(x\). On the right, each series is summed at
\(x\) and \(f_l\) is evaluated at
\(\operatorname{ev}_x(a_0),\dots,\operatorname{ev}_x(a_{L-1})\).
The upper triangle commutes because \(\operatorname{ev}_a(X_i)=a_i\).
For the lower triangle, the two \(\R\)-algebra homomorphisms from
\(\R[X]\) to \(Y\) agree on each variable, with value
\(\operatorname{ev}_x(a_i)\). Their equality
is~\eqref{eq:aeval-natural} with \(h=\operatorname{ev}_x\).
The growth condition on the powers of \(x\) and the completeness of \(Y\)
give the bounded evaluation map. Once that map is given, the substitution
identity is algebraic. For \(Y=\R\) and a real point \(t\) with
\(|t|\le\nu\), the right vertical arrow is \(\lean{eval}(\cdot,t)\) of
Theorem~\ref{thm:mertens}.

\paragraph{Formal coverage.}
The library proves the bounds, atom formulas and uniqueness of the lifts
in Section~\ref{sec:algebras}. The generic algebra lift assumes
\lean{lpOneAlgConvCompat}. Constructing a lift requires a complete
codomain. Proving two existing maps equal by their values on the atoms,
or on \(\delta_1\) in the geometric case, does not. The converse formulas
\lean{eq\_liftCLM\_of\_cols}, \lean{eq\_liftAlgHom\_of\_atoms} and
\lean{eq\_geomAeval\_gen} identify a given map with a constructed lift
and again assume completeness.
For \(\|x\|\le\nu\), \lean{aeval\_of\_norm\_le} uses
\(\|1_Y\|=1\) (\lean{NormOneClass}) to obtain the growth bound with
\(C=1\). The corresponding estimate is
\(\|\operatorname{ev}_x(a)\|\le\|a\|\).

The coproduct, completion, adjunctions and natural bijections above
describe these concrete maps. The project does not package its weighted
lifts as categorical adjunctions or equivalences of hom sets.
The lift norm formulas stated as suprema in Section~\ref{sec:algebras}
follow from the upper bounds and evaluation on the atoms.

\section{Weighted \texorpdfstring{\(\ell^p\)}{lp} spaces}\label{app:lp}

The Taylor examples in Section~\ref{sec:examples} use the Cauchy product
on weighted \(\ell^1\). The Fourier examples in Section~\ref{sec:discussion}
use weighted \(\ell^2\). For weighted \(\ell^p\), two questions arise.
When does convolution give a bounded bilinear product on the space?
How can a certificate bound the norm of a linear operator on it?
We indicate which results the library proves at revision \lean{\LibCommit}.

\paragraph{The space.}
Let \(M\) be a set, \(m\colon M\to(0,\infty)\) a weight and
\(1\le p\le\infty\). The space \(\ell^p_m(M)\) consists of the real or
complex scalar families \(x=(x_j)_{j\in M}\) with
\[
 \|x\|_{p,m}=\Bigl(\sum_{j\in M}\bigl(|x_j|\,m_j\bigr)^p\Bigr)^{1/p}<\infty\quad(p<\infty),\qquad
 \|x\|_{\infty,m}=\sup_{j\in M}|x_j|\,m_j<\infty .
\]
It is a Banach space, the atom \(\delta_j\) has norm \(m_j\), and for
\(p<\infty\) the finitely supported families are dense. At \(p=1\) it is the
general weighted form of the space in Section~\ref{sec:weighted-space}.
In the library, \lean{WeightedLp} is an abbreviation for Mathlib's \lean{lp}.
\begin{lstlisting}
abbrev WeightedLp (M : Type*) (E : M → Type*)
    [∀ m, NormedAddCommGroup (E m)] (p : ℝ≥0∞) : Type _ :=
  lp E p
\end{lstlisting}
Here the fibers \(E(j)\) are copies of the scalars with the weighted norms
and identifications \(\theta_j\) of Section~\ref{sec:weighted-space}.
They are complete, so Mathlib gives completeness of the space.
For finite \(p\), density is \lean{WeightedLp.}\allowbreak\lean{dense\_span\_single}.
The definition specifies no multiplication.

\paragraph{When convolution is a product.}
Let \(M\) be an additive monoid, \(1<p<\infty\) and \(1/p+1/p'=1\).
Suppose the ratios \(m_n/(m_i m_j)\) for \(i+j=n\) satisfy the following
uniform bound. For each \(n\) the sum converges, and for some \(K\ge0\)
independent of \(n\),
\[
 R_{m,p}(n)=\sum_{i+j=n}\Bigl(\frac{m_n}{m_i\,m_j}\Bigr)^{p'}\le K^{p'}
 \qquad(n\in M).
\]
Then for \(x,y\in\ell^p_m(M)\) the convolution
\((x\star y)_n=\sum_{i+j=n}x_iy_j\) converges absolutely at every \(n\),
lies in \(\ell^p_m(M)\), and
\[
 \|x\star y\|_{p,m}\le K\,\|x\|_{p,m}\,\|y\|_{p,m}.
\]
Apply H\"older's inequality on the set \(\{(i,j):i+j=n\}\).
With \(u_i=|x_i|m_i\) and \(v_j=|y_j|m_j\),
\(m_n\sum_{i+j=n}|x_i||y_j|=\sum_{i+j=n}\frac{m_n}{m_im_j}u_iv_j
 \le R_{m,p}(n)^{1/p'}\bigl(\sum_{i+j=n}u_i^pv_j^p\bigr)^{1/p}\),
and these sets of pairs partition \(M\times M\). Since
\(R_{m,p}(n)\le K^{p'}\), summing the \(p\)-th powers gives the stated
norm bound. The product is the unique continuous bilinear extension of
finite convolution. It is associative, has unit \(\delta_0\) of norm
\(m_0\), and is commutative when \(M\) is commutative. The proof requires
neither submultiplicativity of \(m\) nor summability of \(1/m\).

For \(p=2\), suppose \(m_0=1\) and choose \(g\ne0\) in \(M\). For \(t>0\)
with \(z=t^2m_g^2<2\), put \(x=\delta_0+t\delta_g\).
The entries of \(x\star x\) have modulus at least \(1\) at \(0\) and
\(2t\) at \(g\). Thus
\[
 \|x\star x\|_{2,m}^2\ge1+4z>(1+z)^2=\|x\|_{2,m}^4.
\]
No weight with \(m_0=1\) makes convolution contractive on this space.
When the uniform ratio bound holds, \(K\) bounds the bilinear norm.

On \(\Z\), the weights \(m_n=(1+|n|)^\alpha\) satisfy the ratio bound for
\(p=2\) when \(\alpha>1/2\). The triangle inequality and convexity of
\(t\mapsto t^{2\alpha}\) give
\[
 R_{m,2}(n)\le 2^{2\alpha}\sum_{k\in\Z}(1+|k|)^{-2\alpha}<\infty .
\]
Under the uniform ratio bound, the library constructs the bounded product
\lean{WeightedLp.}\allowbreak\lean{convCLM} and proves its uniqueness
and algebraic laws. A Duffing certificate uses this product for
\(m_n=(1+|n|)^2\).

The geometric weights \(m_n=\nu^{|n|}\) for \(\nu\ge1\) restrict on
\(\N\) to the weight \(\nu^n\) of \(\ell^1_\nu\). For \(n\ge0\), the ratio
\(m_n/(m_{n-k}m_k)\) equals \(1\) for \(0\le k\le n\). Thus
\(R_{m,p}(n)\ge n+1\), so no uniform ratio bound holds for
\(1<p<\infty\). Failure of this sufficient condition alone does not
rule out a bounded product. At \(p=2\), for each integer \(N\ge1\) let
\(a_N\) have coordinates \(\nu^{-k}/\sqrt N\) for \(1\le k\le N\) and
zero elsewhere. Then
\[
 \|a_N\|_{2,m}=1,\qquad
 \|a_N\star a_N\|_{2,m}^2=\tfrac{2N}{3}+\tfrac{1}{3N},
\]
so convolution is unbounded in the weighted \(\ell^2\) norm.
For this weight, \lean{WeightedLp.}\allowbreak\lean{Bilateral.prodCLM}
maps the two inputs to \(\ell^\infty_m(\Z)\) with bound \(1\). The Fourier
construction composes this product with a bounded linear map from
\(\ell^\infty_m(\Z)\) back to \(\ell^2_m(\Z)\).

\paragraph{The action of the coefficient algebra.}
Let \(M\) be an additive commutative group and put \(X_p=\ell^p_m(M)\).
Let \(\omega\) be a positive submultiplicative weight, so that
\(\ell^1_\omega(M)\) is the coefficient algebra of
Appendix~\ref{app:categorical}. Assume that \(m\) is \emph{moderate}
with respect to \(\omega\), meaning that for some \(C\ge0\),
\[
 m_{j+l}\le C\,\omega_l\,m_j\qquad(j,l\in M).
\]
The shift \((\tau_lx)_j=x_{j-l}\) then satisfies
\(\|\tau_l\|\le C\omega_l\), with \(\tau_0=I\) and
\(\tau_{j+l}=\tau_j\tau_l\).
The algebra lift applied to the family of shifts gives a bounded unital
homomorphism \(\rho\colon\ell^1_\omega(M)\to\operatorname{End}(X_p)\) into
the algebra of bounded operators on \(X_p\),
\[
 \rho(a)x=\sum_{l\in M}a_l\,\tau_lx,\qquad
 \rho(\delta_l)=\tau_l,\qquad \rho(a\star b)=\rho(a)\,\rho(b),\qquad
 \|\rho\|\le C,
\]
for \(1\le p\le\infty\). Even at \(p=\infty\), the images of the atoms
determine \(\rho\) because its domain is weighted \(\ell^1\).
Let \(\mathbf{wMon}_\infty\) have weighted monoids as objects.
Its arrows are monoid homomorphisms \(f\) satisfying
\(w'(f(j))\le C_f w(j)\) for some \(C_f\ge0\), where \(w\) and \(w'\)
are the two weights. The functor
\(U_\infty\colon\mathbf{BanAlg}_\infty\to\mathbf{wMon}_\infty\) keeps
the multiplication and the norm. In the left triangle the shifts give
the arrow \(\tau\) with bound \(C\). On the right, \(\rho\) is a bounded
algebra homomorphism.
\[
\begin{tikzcd}[column sep=large,row sep=large,ampersand replacement=\&]
 (\operatorname{Mult}(M),\omega)\arrow[r,"\eta"]\arrow[dr,"\tau"']
 \&U_\infty(\ell^1_\omega(M))\arrow[d,"U_\infty(\rho)",dashed]
 \&\ell^1_\omega(M)\arrow[d,"{\exists!\,\rho}",dashed]\\
 \&U_\infty(\operatorname{End}(X_p))\&\operatorname{End}(X_p)
\end{tikzcd}
\]
The arrow \(\rho\) is \lean{WeightedLp.}\allowbreak\lean{symbolAction},
defined by applying \lean{lpOneAlg.liftAlgHom} to the shifts. The shift
bound is \lean{WeightedLp.}\allowbreak\lean{norm\_shift\_le}. The lemma
\lean{WeightedLp.}\allowbreak\lean{symbolAction\(_0\)\_eq} identifies the
lift with the series above. The bound on \(\rho\) is the norm bound on
its underlying linear map. The library constructs this action over a
complete nontrivially normed field.

When \(m\) is submultiplicative, one may take \(\omega=m\) and \(C=1\).
For \(1<p<\infty\), suppose \(m\) also satisfies the uniform ratio bound.
Identifying \(b\in\ell^1_m(M)\) with the element of \(X_p\) having the
same coefficients, one then has \(\rho(b)x=b\star x\)
(\lean{WeightedLp.}\allowbreak\lean{convCLM\_eq\_symbolAction\(_0\)}).

\paragraph{Operators out of \(X_p\).}
Take real scalars and a real Banach space \(Y\).
At \(p=1\), a family \((y_j)_{j\in M}\) gives the columns of a bounded
linear map \(\ell^1_m(M)\to Y\) precisely when the growth condition
\(\|y_j\|\le C m_j\) holds for some \(C\ge0\).
Section~\ref{sec:taylor-services} uses this column bound to estimate the
scalar Taylor tail operator. For \(1<p<\infty\), the atoms still span a
dense subspace of \(X_p\), but their individual norms do not suffice.
The family gives the columns of a bounded linear map precisely when
it satisfies the following \emph{finite synthesis test} for some \(C\ge0\),
\[
 \Bigl\|\sum_{j\in s}c_j\,y_j\Bigr\|\le C\,\Bigl\|\sum_{j\in s}c_j\,\delta_j\Bigr\|_{p,m}
 \qquad(s\subset M\text{ finite},\ c\colon M\to\R),
\]
in which case exactly one bounded linear map \(T\colon X_p\to Y\) has
\(T\delta_j=y_j\), and \(\|T\|\) is the least such \(C\).
The library constructs this map as the extension
\begin{lstlisting}
def liftCLM (hp : p ≠ ∞) (h : LiftBound E p y C) :
    WeightedLp J E p →L[��] Y
\end{lstlisting}
Here \(J=M\), \(\Bbbk=\R\), \(Y\) is complete, and \lean{LiftBound}
records the displayed finite synthesis test. The theorem
\lean{WeightedLp.}\allowbreak\lean{eq\_liftCLM} gives uniqueness, and
\lean{WeightedLp.}\allowbreak\lean{opNorm\_le\_iff\_liftBound}
characterizes the operator norm by the same test.
At \(p=1\) the finite synthesis test reduces to the column bound.
For \(1<p<\infty\) take an infinite \(M\), \(m\equiv1\) and \(y_j=u\)
with \(\|u\|=1\). Every column has norm \(1\), but taking \(c_j=1\) on
\(N\) indices gives \(N\) on the left of the test and \(C N^{1/p}\)
on the right. No fixed \(C\) works for all \(N\).
The test characterizes the bounded operators. On an infinite index set,
neither testing columns separately nor testing all finite combinations
is by itself a finite computation. Particular operators admit simpler
bounds.

For bounded diagonal entries \(d_j\), the operator norm is exactly
\(\sup_j|d_j|\) (\lean{WeightedLp.}\allowbreak\lean{norm\_diagCLM}).
For a matrix \(G=(G_{ij})\) on \(X_p\), let \(R_{\mathrm{row}}\) and
\(R_{\mathrm{col}}\) bound every row sum and every column sum, respectively,
of \(|G_{ij}|\,m_i/m_j\). Then
\(R_{\mathrm{row}}^{1/p'}R_{\mathrm{col}}^{1/p}\) bounds the operator norm
(\lean{WeightedLp.}\allowbreak\lean{schur\_bound\_of}). For
\(1<p<\infty\) these sum bounds are sufficient but not necessary.
For an operator supported on a finite matrix block at \(p=2\),
\lean{WeightedLp.}\allowbreak\lean{norm\_le\_sqrt\_of\_block} gives the
bound \(\sqrt{R_{\mathrm{row}} R_{\mathrm{col}}}\).

A different estimate applies in \(\ell^2_m\) when two operators have
disjoint output supports. If their norms are at most \(b_1\) and \(b_2\),
the norm of their sum is at most \(\sqrt{b_1^2+b_2^2}\), since the
squared output norms add for each input. The library applies this
estimate to the two parts of the tail defect in
\lean{WeightedLp.}\allowbreak\lean{norm\_tailDefect\_le}.
The counterexamples above, including the family \(a_N\), are given here
as mathematical arguments.

\bibliographystyle{alpha}
\bibliography{references}

\end{document}